\documentclass[11pt,reqno]{amsart}
\usepackage{amssymb,amscd,amsbsy}
\usepackage{amssymb,amscd,amsbsy,mathrsfs}
\usepackage{mathdots}
\newcommand{\lb}{\linebreak}

\renewcommand{\a}{\alpha}
\renewcommand{\b}{\beta}

\renewcommand{\d}{\delta}
\newcommand{\e}{\varepsilon}

\newcommand{\z}{\zeta}

\renewcommand{\l}{\lambda}

\newcommand{\s}{\sigma}
\renewcommand{\t}{\tau}

\newcommand{\f}{\varphi}
\renewcommand{\o}{\omega}

\newcommand{\D}{\Delta}
\renewcommand{\L}{\Lambda}

\newcommand{\E}{{\mathscr E}}

\newcommand{\h}{{\mathscr H}}
\newcommand{\I}{{\mathscr I}}

\newcommand{\K}{{\mathscr K}}
\newcommand{\cL}{{\mathscr L}}

\newcommand{\X}{{\mathscr X}}
\newcommand{\Y}{{\mathscr Y}}

\newcommand{\C}{{\Bbb C}}
\newcommand{\T}{{\Bbb T}}
\newcommand{\pp}{{\Bbb P}}
\newcommand{\dd}{{\Bbb D}}
\newcommand{\R}{{\Bbb R}}
\newcommand{\Z}{{\Bbb Z}}

\newcommand{\0}{{\boldsymbol{0}}}

\newcommand{\bs}{\boldsymbol}

\newcommand{\m}{{\boldsymbol m}}
\newcommand{\bS}{{\boldsymbol S}}

\newcommand{\rf}[1]{(\ref{#1})}

\newcommand{\df}{\stackrel{\mathrm{def}}{=}}

\newcommand{\Ker}{\operatorname{Ker}}
\newcommand{\re}{\operatorname{Re}}
\newcommand{\spn}{\operatorname{span}}

\newcommand{\clos}{\operatorname{clos}}

\newcommand{\trace}{\operatorname{trace}}

\newcommand{\const}{\operatorname{const}}

\newcommand{\eeq}{\end{equation}}
\newcommand{\beq}{\begin{equation}}
\newcommand{\bay}{\begin{eqnarray}}
\newcommand{\ba}{\begin{align*}}
\newcommand{\ea}{\end{align*}}
\newcommand{\ey}{\end{eqnarray}}
\newcommand{\bey}{\begin{eqnarray*}}
\newcommand{\eey}{\end{eqnarray*}}

\newcommand{\be}{\infty}

\newcommand{\bl}{\blacksquare}

\newcommand{\Range}{\operatorname{Range}}

\newcommand{\Pf}{{\bf Proof. }}
\newcommand{\im}{\operatorname{Im}}
\renewcommand{\re}{\operatorname{Re}}
\newcommand{\ov}{\overline}

\newtheorem{thm}{\hspace{\parindent}Theorem}[section]

\newtheorem{cor}[thm]{\hspace{\parindent}Corollary}
\newtheorem{lem}[thm]{\hspace{\parindent}Lemma}

\theoremstyle{remark}

\newtheorem*{rem*}{Remark}

\newcommand\CA{{\rm C}_{\rm A}}
\newcommand{\OL}{{\rm OL}}
\newcommand{\OLr}{{\rm OL}^{\rm res}}
\newcommand{\OLA}{{\rm OL}_{\rm A}}
\newcommand{\OLAr}{{\rm OL}_{\rm A}^{\rm res}}

\newcommand{\ri}{{\rm i}}

\newcommand\mM{\mathcal{M}}

\newcommand\dg{\frak D}

\begin{document}

\newcommand{\vse}{\vspace{.2in}}
\numberwithin{equation}{section}

\title{Absolute continuity of spectral shift}
\thanks{The publication was prepared with the support of the
"RUDN University Program 5-100"}

\author{M.M. Malamud, H. Neidhardt and V.V. Peller}
%\thanks{the author is partially supported by NSF grant DMS 1300924}
%\thanks{Corresponding author: V.V. Peller; email: peller@math.msu.edu}

\maketitle

\begin{abstract}
In this paper we develop the method of double operator integrals to prove trace formulae for functions of contractions, dissipative operators, unitary operators and 
self-adjoint operators. To establish the absolute continuity of spectral shift, we use the Sz.-Nagy theorem on the absolute continuity of the spectrum of the minimal unitary dilation of a completely nonunitary contraction.
We also give a construction of an intermediate contraction for a pair of contractions with trace class difference.
\end{abstract}

\maketitle

\tableofcontents
\normalsize
\vspace*{-1cm}

\setcounter{section}{0}
\section{\bf Introduction}
\setcounter{equation}{0}
\label{In}

\

The notion of spectral shift function was introduced by physicist I.M. Lifshits in \cite{L}. Later M.G. Krein elaborated the notion of spectral shift function in \cite{Kr} a most general situation; he showed that for self-adjoint operators $A_0$ and $A_1$ with trace class difference, there exists a unique real function 
$\bs{\xi}=\bs{\xi}_{A_0,A_1}$ in $L^1(\R)$ (it is called the {\it spectral shift function for the pair} $\{A_0,A_1\}$) such that the following trace formula holds:
\bay
\label{foslsso}
\trace\big(f(A_1)-f(A_0)\big)=\int_\R f'(t)\bs{\xi}_{A_0,A_1}(t)\,dt
\ey
%for all functions $f$ whose derivative $f'$ is the Fourier transform of a complex Borel measure. 
for sufficiently nice functions $f$.
To prove the existence of the spectral shift function, he
introduced  the concept of {\it perturbation determinant} $\Delta_{A_1/A_0}$ and proved
the inversion formula
 $$
 \bs{\xi}_{A_0,A_1}(t) = \frac{1}{\pi}\lim_{y\downarrow 0}\im (\log(\Delta_{A_1/A_0}(t +\ri y)))
\quad \mbox{for a.e.} \quad  t \in \R,
$$
where  $\Delta_{A_1/A_0}(\z)\df\det(I + (A_1-A_0)(A_0-\z I)^{-1})$ (see \cite{Kr} and
\cite{Ya}).

Another approach to trace formula was given by Birman and Solomyak in \cite{BS2}. Their approach is based on differentiating in the trace norm $\bS_1$ of the parametric family $f(A_t)-f(A_0)$ and computing the trace of the double operator integrals that represent the derivative. 
Here $A_t\df(1-t)A_0+tA_1$, $0\le t\le1$.
However, the approach of Birman and Solomyak did not lead to the absolute continuity of spectral shift. They showed that
\bay
\label{foslemer}
\trace\big(f(A_1)-f(A_0)\big)=\int_\R f'(t)\,d\nu(t)
\ey
for a finite signed Borel measure $\nu$ on $\R$ that is uniquely determined by the pair $\{A_0,A_1\}$
and for sufficiently nice functions $f$. Such a measure $\nu$ can be called the {\it spectral shift measure}. However, it follows from the Krein theorem that is must be absolutely continuous with respect to Lebesgue measure $\m$ and 
$d\nu=\bs\xi_{A_0,A_1}\,d\m$.

Let us also mention here the paper \cite{PSZ}, in which the authors give another proof of the Lifshits--Krein trace formula.

Note that spectral shift function plays an important role in perturbation theory. We mention here
the paper \cite{BK}, in which the following important formula was found:
$$
\det{\mathscr S}(x)=e^{-2\pi{\rm i}\xi(x)},
$$
where ${\mathscr S}$ is the scattering matrix corresponding to the pair 
$\{A_0,A_1\}$.

Krein extended in \cite{Kr3} formula \eqref{foslsso} to the class ${\mathcal W}_1(\R)$  of functions whose derivative is the Fourier transform of a complex Borel measure.
He also observed in \cite{Kr} that the right-hand side of \rf{foslsso} makes sense for arbitrary Lipschitz functions $f$ and posed the problem to describe the class of functions, for which formula \rf{foslsso} holds for all pairs of self-adjoint operators with difference of trace class $\bS_1$.
It turned out that trace formula \rf{foslsso} cannot be generalized to the class of all Lipschitz functions. Indeed, it was shown in \cite{F} that there exist a Lipschitz function $f$ on $\R$ and self-adjoint operators $A_1$ and $A_0$ such that $A_1-A_0\in\bS_1$, but $f(A_1)-f(A_0)\not\in\bS_1$. In \cite{Pe1} and \cite{Pe3} it was proved that \rf{foslsso} holds for functions $f$ in the (homogeneous) Besov space
$B_{\be,1}^1(\R)$ and does not hold unless $f$ locally belongs to the Besov space
$B_{1,1}^1$

Krein's problem was completely solved recently in \cite{Pe6}. It was shown in \cite{Pe6} that trace formula \rf{foslsso} holds for arbitrary pairs $\{A_0,A_1\}$ of not necessarily bounded self-adjoint operators with trace class difference if and only if $f$ is an {\it operator Lipschitz function}, i.e., the inequality 
$$
\|f(A)-f(B)\|\le\const\|A-B\|
$$
holds for arbitrary self-adjoint operators $A$ and $B$. To solve the Krein problem
in \cite{Pe6}, the method of Birman and Solomyak \cite{BS2} was modified.
In particular, instead of differentiating the function $t\mapsto f(A_t)-f(A_0)$ in the trace norm, 
the derivative is taken in the Hilbert--Schmidt norm $\bS_2$.
This allowed the author to extend the Lifshits--Krein trace formula to the class of all operator Lipschitz functions. 

However, the methods of \cite{Pe6} did not lead to the absolute continuity of spectral shift. To prove the validity of trace formula \rf{foslsso} for all operator Lipschitz functions, Krein's theorem had to be used which implied that the measure $\nu$ in formula \rf{foslemer} is absolutely continuous with respect to Lebesgue measure.

In this paper we give a further development of the Birman--Solomyak idea \cite{BS2} and {\it get directly the absolute continuity of spectral shift}. However, it turned out that to achieve this purpose, we have to consider similar problems not only for self-adjoint operators. To get absolute continuity, we start with the case of functions of contractions. Then we use the result for contractions to proceed to functions of unitary operators. Finally, this will allow us to treat the cases of functions of self-adjoint operators and functions of maximal dissipative operators. 

Recall that Krein introduced in \cite{Kr2} the notion of spectral shift function for pairs of unitary operators with trace class difference. Namely, for a pair of unitary operators $\{U_0,U_1\}$ with trace class difference, he proved that there exists a real function $\bs{\xi}$ in $L^1(\T)$, unique modulo an additive constant, (called a {\it spectral shift function for} $\{U_0,U_1\}$ ) such that the trace formula 
\bay
\label{fosuo}
\trace\big(f(U_1)-f(U_0)\big)=\int_\T f'(\z)\bs{\xi}(\z)\,d\z
\ey
holds for functions $f$ whose derivative $f'$ has absolutely convergent Fourier series.

In the recent paper \cite{AP+} an analog of the result of \cite{Pe6} was obtained: 
the maximal class of functions $f$, for which formula \rf{fosuo} holds for arbitrary pairs $\{U_0,U_1\}$ of unitary operators with $U_1-U_0\in\bS_1$ coincides with the class $\OL(\T)$ of operator Lipschitz functions on $\T$. This class can be defined by analogy with the definition of the class of operator Lipschitz functions on $\R$. The authors of \cite{AP+} considered the parametric family $e^{{\rm i}tA}U_0$, $0\le t\le1$, where $A$ is a trace class self-adjoint operator such that 
$U_1=e^{{\rm i}A}U_0$. However, unlike in the case of functions of self-adjoint operators (see \cite{Pe6}), it is still unknown whether
the function $t\mapsto f(e^{{\rm i}tA})$ is differentiable in the $\bS_2$ norm for an arbitrary operator Lipschitz function $f$ on the unit circle $T$. Instead, it was shown in \cite{AP+} that this functions must be differentiable in the strong operator topology. This allowed the authors to extend trace formula \rf{fosuo} to the class of all operator Lipschitz functions on $\T$. To be more precise, the methods of \cite{AP+} lead to the trace formula
$$
\trace\big(f(U_1)-f(U_0)\big)=\int_\T f'(\z)\,d\nu(\z)
$$
for a finite signed Borel measure $\nu$ on $\T$ and for arbitrary operator Lipschitz
functions $f$ on $\T$. As in the case of functions of self-adjoint operators the absolute continuity of $\nu$ follows from Krein's theorem of \cite{Kr2}.

In the recent paper \cite{MNP2} (see also \cite{MNP1}) definitive results were obtained on trace formulae for functions of contractions.
For a pair of Hilbert space contractions $\{T_0,T_1\}$ with trace class difference $T_1-T_0$, we considered the parametric family $T_t=T_0+(1-t)(T_1-T_0)$, $0\le t\le1$, of contractions.
By differentiating the parametric family $f(T_t)$ and computing the trace of certain double operator integrals with respect to semi-spectral measures we showed in 
\cite{MNP2} that there exists a complex Borel measure $\nu$ on $\T$ such that 
$$
\trace\big(f(T_1)-f(T_0)\big)=\int_\T f'(\z)\,d\nu(\z)
$$
for arbitrary operator Lipschitz functions $f$ analytic in $\dd$. However, we did not 
get the absolute continuity of $\nu$.
To prove that $\nu$ is absolutely continuous, we used earlier results of \cite{MN} for functions of resolvent comparable maximal dissipative operators under certain additional restrictions. Similar results were also obtained in \cite{MNP2} for functions of resolvent comparable maximal dissipative operators.

Note that the space of operator Lipschitz functions is nonseparable. That is why to prove that trace formulae hold for all operator Lipschitz functions, it is not enough to prove them for nice functions.

%The Sz.-Nagy--Foia\c s functional calculus (see \cite{SNF})
%associates with each function $f$ in the {\it disk-algebra} $\CA$ (i.e., the space of functions analytic in the disk $\dd$ and continuous in its closure) the operator $f(T)$. This functional calculus is linear and multiplicative and the von Neumann inequality 
%$
%\|f(T)\|\le\max\{|f(\z)|:~|\z|\le1\}
%$
%holds for $f\in\CA$.

%When we proceed to the case of contractions, it is natural to divide the problem in  two
%problems: (i) the existence of a spectral shift function and a trace formula for resolvents; (ii) a
%description of the largest class  of functions for which the trace formula 
%holds.

%First generalizations  of formula \eqref{foslsso} to the case of accumulative
%(dissipative) operators were obtained by Rybkin \cite{Ryb84} and Krein \cite{Kr87}. 
%First generalizations  of formula \eqref{foslsso} to the case of pairs $\{A_0, A_1\}$ with
%an m-accumul\-ative (dissipative) operator $A_1$ were obtained by Rybkin \cite{Ryb84} and
%Krein \cite{Kr87}.
%For
%instance, Krein treated the case when $A_0 = A^*_0$,  $A_1 = A_0 - {\rm i}V$,
%$V\ge 0$ and $V\in \bS_1$ and proved in \cite{Kr87} an analog of formula \eqref{foslsso} with right-hand side $\int_\R f'(t)\,d\nu(t)$ for
% a complex Borel measure $\nu$ and for functions $f$ of class $W_1^+(\R)$, i.e., functions $f$ whose derivative is the Fourier transform of a complex measure supported in $[0,\be)$.
% $f\in\mathcal K_+$ of $C^1_{loc}$ functions whose derivatives 
% $f'(z)=
%\int_{\R_+}e^{-izt}\;dp(t)$, $z\in \clos\C_+$, are the Laplace transform of finite
%measures. 

In this paper we consistently use differentiation of parametric families of functions of operators and writing explicit formulae for the trace of the derivatives in terms of double operator integrals.
Moreover, we manage
not only to describe the maximal classes of functions, for which, the corresponding trace formula is applicable, but also establish directly the absolute continuity of spectral shift.

We start with the case of functions of contractions. In \S\;\ref{szhali}  we deduce the absolute continuity of spectral shift for contractions from the Sz.-Nagy theorem on the absolute continuity of the spectrum of minimal unitary dilations of completely nonunitary contractions. This will allow us to obtain in \S\;\ref{Uni} the desired result for functions of unitary operators by applying the brothers Riesz theorem. Then we use Cayley transform to obtain in \S\;\ref{samosoop} similar results for functions of self-adjoint operators. Finally, in \S\;\ref{rezsradis} and
\S\;\ref{disaddi} we apply Cayley transform to contractions and obtain definitive results for functions of maximal dissipative operators. 

Note that a spectral shift function for a pair of contractions is not unique. If 
$\bs\xi$ is a spectral shift function for a pair $\{T_0,T_1\}$ of contractions with trace class difference, then all spectral shift functions for $\{T_0,T_1\}$ can be parametrized by
$\bs\xi+h$, where $h$ ranges over the Hardy class $H^1$. It is not always possible to find a real-valued spectral shift function. In \S\;\ref{promzh} for a pair 
$\{T_0,T_1\}$ of contractions with trace class difference under a mild assumption, we construct an intermediate contraction $T$ such that $T-T_0\in\bS_1$ and such that
the pair $\{T_0,T\}$ has a spectral shift function $\bs\xi_0$ with 
$\im\bs\xi_0\ge\0$, while the pair $\{T,T_1\}$ has a spectral shift function 
$\bs{\xi}_1$ with $\im\bs\xi_1\le\0$. Clearly, the function $\bs\xi\df\bs\xi_0+\bs\xi_1$
is a spectral shift function for the pair $\{T_0,T_1\}$. This allows us to establish in \S\:\ref{promzh} that if $U$ is a unitary operator and $T$ is a contraction such that $T-U\in\bS_1$, then the pair $\{U,T\}$ has a spectral shift function $\bs\xi$ with $\im\bs\xi\ge\0$.

In the Appendix we present a useful result on Sch\"affer matrix unitary dilations of contractions. 

Finally, in \S\;\ref{dois} we give an introduction to double operator integrals and in \S\;\ref{OperLirazra} we give an introduction to operator Lipschitz functions.

We consider in this paper only operators on {\it separable} Hilbert spaces.

%The purpose of this paper is to prove that for a pair $\{T_0,T_1\}$  of contractions satisfying $T_0-T_1\in\bS_1$, there exists a function $\bs\xi$ in $L^1(\T)$ (called a {\it spectral shift function for the pair} $\{T_0,T_q\}$ ) such that the following trace formula
%\bay
%\label{tfszha}
%\trace\big(f(T_1)-f(T_0)\big)=\int_\T f'(\z)\bs{\xi}(\z)\,d\z
%\ey
%holds for an arbitrary operator Lipschitz function $f$ analytic in $\dd$.
%
%To obtain the main result, we combine two approaches. The first approach is based on double operator integrals with respect to semi-spectral measures and uses an idea of \cite{BS2}. It leads to a trace formula 
%$\trace\big(f(T)-f(R)\big)=\int_\T f'(\z)\,d\nu(\z)$ for an arbitrary operator Lipschitz function $f$ analytic in $\dd$, where $\nu$ is a Borel measure on $\T$.
%
%The second approach develops ideas of Krein \cite{Kr} - \cite{Kr87} and relies on the study of the perturbation determinant for pairs
%$\{L_0,L_1\}$ of m-accumulative operators and allows us to prove the trace formula for
%resolvents without  additional restrictions on $\{L_0,L_1\}$ imposed in \cite{MalNei2015}.
%
%We denote by $\mB(\h)$ the set of bounded linear operators on a Hilbert space $\h$ and by $\bS_p$ the
%Schatten-von Neumann ideal. In particular,  $\bS_1$ is the trace
%class.  For a closed densely
%defined operator $Q$, we use the notation $\rho(Q)$, $\s_{\rm p}(Q)$,  and 
%$\s_{\rm c}(Q)$
%denote the resolvent set and the point and the continuous spectra  of $Q$. Recall
%that 
%$\s_{\rm c}(Q) =\{\l\notin\s_{\rm p}(Q):~\Range(Q-\l I)\ne\clos\Range(Q-\l I) =
%\h$\}.

\

\section{\bf Double operator integrals and Schur multipliers}
\setcounter{equation}{0}
\label{dois}

\

Double operator integrals
$$
\iint_{\X\times\Y}\Phi(x,y)\,dE_1(x)Q\,dE_2(y)
$$
were introduced by Yu.L. Daletskii and S.G. Krein in \cite{DK}. Later
Birman and Solomyak elaborated their beautiful theory of double operator integrals,
see \cite{BS} and \cite{BS3} (see also \cite{AP} and references therein). Here $\Phi$ is a bounded measurable function, $E_1$ and $E_2$ are spectral measures on Hilbert space
defined on $\s$-algebras of subsets of sets $\X$ and $\Y$ and
$Q$ is a bounded linear operator. 

The starting point of the Birman--Solomyak approach \cite{BS} is the case when $Q\in\bS_2$.
Under this assumption double operator integrals can be defined for arbitrary bounded measurable functions $\Phi$ as follows. Consider the spectral measure ${\mathcal E}$ whose values are orthogonal
projections on the Hilbert space $\bS_2$, which is defined by
$$
{\mathcal E}(\L\times\D)Q=E_1(\L)QE_2(\D),\quad Q\in\bS_2,
$$
where $\L$ and $\D$ are measurable subsets of $\X$ and $\Y$. 
Obviously, left multiplication by $E_1(\L)$ commutes with right multiplication by 
$E_2(\D)$.
It was shown in \cite{BS} (see also \cite{BS4}) that ${\mathcal E}$ extends to a spectral measure on
$\X\times\Y$ if $\Phi$ is a bounded measurable function on $\X\times\Y$ and, by definition,
$$
\iint_{\X\times\Y}\Phi(x,y)\,d E_1(x)Q\,dE_2(y)\df
\left(\,\iint_{\X\times\Y}\Phi\,d{\mathcal E}\right)Q.
$$
Clearly,
$$
\left\|\iint_{\X\times\Y}\Phi(x,y)\,dE_1(x)Q\,dE_2(y)\right\|_{\bS_2}
\le\|\Phi\|_{L^\be}\|Q\|_{\bS_2}.
$$

If $Q$ is an arbitrary bounded linear operator, then for the double operator integral to make sense, $\Phi$ has to be a Schur multiplier with respect to $E_1$ and $E_2$, (see \cite{Pe1} and \cite{AP}). It is well known (see \cite{Pe1}, \cite{AP}, \cite{Pe+} and \cite{Pi}) that $\Phi$ is a Schur multiplier if and only if 
$\Phi$ belongs to the Haagerup tensor product 
$L^\be_{E_1}\otimes_{\rm h}L^\be_{E_2}$ of 
$L^\be_{E_1}$ and $L^\be_{E_2}$ or, in other words, $\Phi$
admits a representation
\bay
\label{htenppre}
\Phi(x,y)=\sum_{n\ge0}\f_n(x)\psi_n(y),
\ey
where the $\f_n$ and $\psi_n$ are measurable functions such that
$$
\sum_{n\ge0}|\f_n|^2\in L^\be_{E_1}\quad\mbox{and}\quad
\sum_{n\ge0}|\psi_n|^2\in L^\be_{E_2}.
$$
In this case
$$
\iint\Phi(x,y)\,dE_1(x)Q\,dE_2(y)=
\sum_{n\ge0}\left(\int\f_n(x)\,dE_1(x)\right)Q\left(\int\psi_n(y)\,dE_2(y)\right),
$$
the series on the right converges in the weak operator topology
and the right-hand side does not depend on the choice of a representation in 
\rf{htenppre}.

In this paper we also need double operator integrals with respect to {\it semi-spectral measures}
\bay
\label{dvoopipolu}
\iint\Phi(x,y)\,d\E_1(x)Q\,d\E_2(y).
\ey
Such double operator integrals were introduced in \cite{Pe3} (see also \cite{Pe6}).
By analogy with the case of double operator integrals with respect to spectral measures, double operator integrals of the form \rf{dvoopipolu} can be defined for arbitrary bounded measurable functions $\Phi$ in the case when $Q\in\bS_2$ and for 
functions $\Phi$ in $L^\be_{\E_1}\otimes_{\rm h}L^\be_{\E_2}$ in the case of an arbitrary bounded operator $Q$.

We refer the reader to the recent surveys \cite{AP} and \cite{Pe+} for detailed information.

Each contraction $T$ (i.e., an operator of norm at most 1) on a Hilbert space $\h$ has a {\it minimal unitary dilation} $U$, i.e., $U$ is a unitary operator on a Hilbert space $\K$, $\K\supset\h$, $T^n=P_\h U^n\big|\h$ for $n\ge0$ and $\K$ is the closed linear span of
$U^n\h$, $n\in\Z$ (see \cite{SNF}, Ch.~I, Th.~4.2). Here $P_\h$ is the orthogonal projection onto $\h$.
The {\it semi-spectral measure $\E_T$ of} $T$ is defined by
\bay
\label{pospecomspe}
\E_T(\D)\df P_\h E_U(\D)\big|\h,
\ey
where $E_U$ is the spectral measure of $U$ and $\D$ is a Borel subset of $\T$. It is easy to see that 
$T^n=\int_\T\z^n\,d\E_T(\z),~ n\ge0.$

Note that if $U$ is a unitary dilation of $T$ on a Hilbert space $\K$, $\K\supset\h$,
that is not necessarily minimal, then formula \rf{pospecomspe} still holds. Indeed,
the closed linear span of $U^n\h$, $n\in\Z$, is a reducing subspace of $U$ that contains $\h$. Thus, for every Borel subset $\D$,
$$
P_\h E_U(\D)\big|\h=P_\h E_{U_0}(\D)\big|\h,
$$
where $U_0$ is the restriction of $U$ to the closed linear span of $U^n\h$, $n\in\Z$.

Note also that for a minimal unitary dilation $U$, the operator measures $\E_T$  and  $E_U$  are spectrally equivalent, i.e.,  $\E_T(\D) = 0$  if and only if $E_U(\D)=0$,  and the multiplicity functions
$N_{\E_T}$ and $N_{E_U}$ coincide almost everywhere with respect to ${E_U}$
(see \cite{MM}).

For a maximal dissipative operator $L$ in a Hilbert space $\h$, its semi-spectral measure $\E_L$ can be defined in the following way. It is well known 
%(see \cite{AP1}) 
that $L$ has a minimal resolvent self-adjoint dilation $A$, i.e., $A$ is a self-adjoint operator in a Hilbert space $\K$, $\h\subset\K$, 
$$
(L-\z I)^{-1}=P_\h(A-\z I)^{-1}\big|\h,\quad\im\z<0,
$$
and $\K=\clos\spn\big\{(A-\z I)^{-1}\h:~\im\z<0\big\}$.

The {\it semi-spectral measure} $\E_L$ of $L$ is defined by
$$
\E_L(\D)\df P_\h E_A(\D)\big|\h
$$
for a Borel subset $\D$ of $\R$, where $E_A$ stands for the spectral measure of $A$.

The functional calculus $f\mapsto f(L)$ can be defined on the class of functions $f$ bounded and analytic in $\C_+$ and such that $f\big|\R$ is continuous on $\R$. We put
$$
f(L)\df\int_\R f(t)\,d\E_L(t)
$$
(see, e.g., \cite{AP1}).

\

\section{\bf Operator Lipschitz functions and divided differences}
\setcounter{equation}{0}
\label{OperLirazra}

\

It is well known (see \cite{AP}) that for a continuous function $f$ on $\R$ the following are equivalent:

\medskip

(i) {\it $f$ is operator Lipschitz, i.e.,
\bay
\label{Lipogr}
\|f(A)-f(B)\|\le\const\|A-B\|
\ey
for arbitrary bounded self-adjoint operators $A$ and $B$};

(ii) {\it inequality {\em\rf{Lipogr}} holds for arbitrary not necessarily bounded self-adjoint operators $A$ and $B$ with bounded $A-B$};

(iii) {\it for arbitrary not necessarily bounded self-adjoint operators $A$ and $B$
such that $A-B\in\bS_1$, the operator $f(A)-f(B)$ is also of trace class}; 

(iv) {\it the inequality
$$
\|f(A)-f(B)\|_{\bS_1}\le\const\|A-B\|_{\bS_1}
$$
holds for arbitrary self-adjoint operators $A$ and $B$ with trace class difference}.

\medskip

We use the notation $\OL(\R)$ for the class of operator Lipschitz functions on $\R$.

It is also well known (see the original paper \cite{JW} and the recent survey \cite{AP}) that operator Lipschitz functions are differentiable everywhere on $\R$.
Moreover, they are also differentiable at infinity, i.e., the limit
$$
\lim_{|t|\to\be}\frac{f(t)}t
$$
exists for each function $f$ in $\OL(\R)$.
However, operator Lipschitz functions are not necessarily continuously differentiable, see \cite{KS} and the survey \cite{AP}.

Similar properties hold for functions of unitary operators. For a continuous function $f$ on the unit circle $\T$ the following are equivalent:

\medskip

(i) {\it $f$ is operator Lipschitz, i.e.,
$$
\|f(U)-f(V)\|\le\const\|U-V\|
$$
for arbitrary unitary operators $U$ and $V$};

(ii) {\it for arbitrary unitary operators $U$ and $V$
such that $U-B\in\bS_1$, the operator $f(U)-f(V)$ is also of trace class}; 

(iii) {\it the inequality
$$
\|f(U)-f(V)\|_{\bS_1}\le\const\|U-V\|_{\bS_1}
$$
holds for arbitrary unitary operators $U$ and $V$ with trace class difference}.

\medskip

This can be proved in a way similar to the case of self-adjoint operators, see \cite{AP}.

It became clear from the papers \cite{DK} and \cite{BS3} that an important role in studying functions of operators under perturbation is played by divided differences. For a differentiable function $f$ on $\R$ the {\it divided difference} $\dg f$ is defined on $\R\times\R$ by
$$
(\dg f)(x,y)\df
\left\{\begin{array}{ll}\displaystyle{\frac{f(x)-f(y)}{x-y}},&x\ne y,\\[.4cm]
f'(x),&x=y.
\end{array}\right.
$$

It is well known (see \cite{AP}) that a differentiable function $f$ on $\R$ is operator Lipschitz if and only if the divided difference $\dg f$ is a Schur multiplier with respect to arbitrary Borel spectral measures on $\R$. Moreover, for an operator Lipschitz function $f$ on $\R$ and for self-adjoint operators $A$ and $B$ with bounded difference the following formula holds:
$$
f(A)-f(B)=\iint_{\R\times\R}(\dg f)(x,y)\,dE_A(x)(A-B)\,dE_B(y),
$$
where $E_A$ and $E_B$ are the spectral measures of $A$ and $B$, see \cite{BS3} and \cite{AP}.

Finally, let us mention that a differentiable function $f$ on $\R$ is operator Lipschitz if and only if there exist sequences $\{\f_n\}_{n\ge1}$ and 
$\{\psi_n\}_{n\ge1}$ of continuous functions on $\R$ such that the limits
$$
\lim_{|x|\to\be}\f_n(x)\quad\mbox{and}\quad\lim_{|x|\to\be}\psi_n(x)
$$
exist,
$$
\sup_{x\in\R}\sum|\f_n(x)|^2<\be
\quad\mbox{and}\quad\sup_{x\in\R}\sum|\psi_n(x)|^2<\be
$$
and
$$
(\dg f)(x,y)=\sum\f_n(x)\psi_n(y),\quad x,~y\in\R,
$$
see \cite{AP}.

Note that similar definitions can be given for functions on the unit circle $\T$ and similar result hold, see \cite{AP}.

Let us proceed now to operator Lipschitz functions on the unit disk and on the upper half-plane.

A function $f$ defined on the closed unit disk $\clos\dd$ is called {\it operator Lipschitz} if 
$$
\|f(N_1)-f(N_2)\|\le\const\|N_1-N_2\|
$$
for arbitrary normal operators $N_1$ and $N_2$ with spectra in $\clos\dd$. 
If $f$ is a function in the disk-algebra $\CA$, then it is operator Lipschitz if and only if 
$$
\|f(T)-f(R)\|\le\const\|T-R\|
$$
for arbitrary contractions $T$ and $R$, see \cite{KS2} and \cite{AP}. 
We use the notation 
$\OLA$ for the class of operator Lipschitz functions in $\CA$. It was proved
in \cite{KS2} that $\OLA=\OL(\T)\cap\CA$, see also\cite{AP}, \S\:3.9.

Note also that for $f\in\OLA$ and for contractions $T$ and $R$ the condition $T-R\in\bS_1$ implies $f(T)-f(R)\in\bS_1$ (see \cite{AP}).

For $f\in\OLA$, we consider the divided difference $\dg f$ defined on 
$\clos\dd\times\clos\dd$ by
\bay
\label{razrakr}
(\dg f)(\z,\t)\df
\left\{\begin{array}{ll}\displaystyle{\frac{f(\z)-f(\t)}{\z-\t}},&\z\ne \t,\\[.4cm]
f'(\z),&\z=\t.
\end{array}\right.
\ey

We are going to use the following characterization of the divided difference 
$\dg f$ for functions in $\OLA$ (see \cite{AP}, Theorems 3.9.1 and 3.9.2):

\medskip

{\it Let $f$ be a function analytic in $\dd$. Then $f\in\OLA$ if and only if
$\dg f$ admits a representation 
\bay
\label{Haagerpre}
(\dg f)(z,w)=\sum_{n\ge1}\f_n(z)\psi_n(w),\quad z,\;w\in\clos\dd,
\ey
where $\f_n$ and $\psi_n$ are functions in $\CA$  such that
\bay
\label{proisumm}
\Big(\sup_{z\in\dd}\sum_{n\ge1}|\f_n(z)|^2\Big)\Big(\sup_{w\in\dd}\sum_{n\ge1}|\psi_n(w)|^2\Big)<\be.
\ey
If $f\in\OLA$, then the functions $\f_n$ and $\psi_n$ can be chosen so that
the left-hand side of {\em\rf{proisumm}} is equal to $\|f\|_{\OLA}$.}

\medskip

Let $T_0$ and $T_1$ be contractions on Hilbert space, and let $\E_0$ and $\E_1$ be their semi-spectral measures. Suppose now that $f\in\OLA$. 
Consider a representation of $\dg f$ in the form \rf{Haagerpre}, where $\f_n$ and $\psi_n$ are functions in $\CA$ satisfying \rf{proisumm}. Then for a bounded linear operator $K$, we have
\bay
\label{dvoopiryad}
\iint_{\T\times\T}\big(\dg f)(\z,\t)\,d\E_1(\z)K\,d\E_0(\t)
=\sum_{n=1}^\be\f_n(T_1)K\psi_n(T_0)
\ey
(see Section 3.9 of \cite{AP}).
This implies (see Theorem 3.9.9 of \cite{AP}) that
\bay
\label{razdvoopi}
f(T_1)-f(T_0)=
\iint_{\T\times\T}\big(\dg f)(\z,\t)\,d\E_1(\z)(T_1-T_0)\,d\E_0(\t).
\ey

Similarly, a function $f$ on the closed upper half-plane $\clos\C_+$ is called
operator Lipschitz if the inequality
$$
\|f(N_1)-f(N_2)\|\le\const\|N_1-N_2\|
$$
holds for arbitrary normal operators $N_1$ and $N_2$ with spectra in $\clos\C_+$.
Consider the class of operator Lipschitz functions on $\clos\C_+$ that are analytic in $\C_+$ and denote it by
$\OLA(\C_+)$. 

It is also well known (\cite{JW}, see also \cite{AP}) that for an arbitrary function $f$ in $\OLA(\C_+)$, the limit
\bay
\label{pronabe}
\lim_{\im\z\ge0,\,|\z|\to\be}\frac{f(\z)}{\z}
\ey
exists.

We use the notation $\CA(\C_+)$ for the set of functions analytic in $\C_+$, continuous in $\clos\C_+$ and having finite limit at infinity.
For a function $f$ in $\OLA(\C_+)$, we can define the divided difference $\dg f$
as in \rf{razrakr}. The following analog of the above result holds (see \cite{AP}, Th. 3.9.6):

\medskip

{\it Let $f\in\OLA(\C_+)$. Then there are sequences $\{\f_n\}_{n=1}^\be$ and 
$\{\psi_n\}_{n=1}^\be$ in $\CA(\C_+)$ such that$$
\left(\sup_{z\in\C_+}\sum_{n=1}^\be|\f_n(z)|^2\right)\left(\sup_{w\in\C_+}\sum_{n=1}^\be|\psi_n(w)|^2\right)
=\|f\|_{\OL(\C_+)}^2
$$
and
$$
(\dg f)(z,w)=\sum_{n=1}^\be\f_n(z)\psi_n(w).
$$
Herewith the series converge uniformly while $z$ and $w$ range over compact subsets of the open upper half-plane}.

\

\section{\bf An analog of the Lifshits--Krein trace formula for contractions}
\setcounter{equation}{0}
\label{szhali}

\

Recall that in our paper \cite{MNP2} the method of \cite{Pe6} (which is in turn based on an idea of \cite{BS3}) was extended to the case of functions of contractions which led to the following trace formula
\bay
\label{maksborme}
\trace\big(f(T_1)-f(T_0)\big)=\int_\T f'(\z)\,d\nu(\z)
\ey
for arbitrary operator Lipschitz functions $f$ in $\OLA$ and for a complex Borel measure $\nu$ on $\T$. Then another, more complicated, method was used in \cite{MNP2} to prove that the measure $\nu$ must be absolutely continuous with respect to normalized Lebesgue measure $\m$ on $\T$. That method is based on a development of Krein's method of perturbation determinants \cite{Kr}, \cite{Kr2} and \cite{Kr3}.

In this section we show that the method of differentiating parametric families of contractions and computing the trace of the corresponding double operator integrals
can give more. Not only allows it us to describe the maximal class of functions $f$, for which the trace formula is applicable, but also it can lead to the absolute continuity of the measure $\nu$.
We also obtain at the end of the section new trace formulae in terms of A-integrals.

Later we will see that the absolute continuity of spectral shift in the case of functions of contractions can allow us to establish absolute continuity for self-adjoint operators, unitary operators and dissipative operators.

However, it is important to start with the case of functions of contractions and deduce from this case the results for the other classes of operators.

We are going to state here certain results obtained in \cite{MNP2} and to repeat the reasoning that was used in \cite{MNP2} to obtain formula \rf{maksborme} for functions  $f$ in $\OLA$.

The following result is Theorem 4.1 of \cite{MNP2}:

\medskip

{\it Let $f\in\OLA$ and let $T_0$ and $T_1$ be contractions on Hilbert space and
$T_t\df T_0+\lb t(T_1-T_0)$, $0\le t\le1$. Then
$$
\lim_{s\to0}\frac1s\big(f(T_{t+s})-f(T_t)\big)
=\iint_{\T\times\T}\big(\dg f)(\z,\t)\,d\E_t(\z)(T_1-T_0)\,d\E_t(\t)
$$
in the strong operator topology,
where $\E_t$ is the semi-spectral measure of $T_t$}.

\medskip

Note that we are going to prove an analog of Theorem 4.1 of \cite{MNP2} for dissipative operators in \S\;8 of this paper.

We also need Theorem 5.1 of \cite{MNP2}:

\medskip

{\it Let $f\in\OLA$. Suppose that $T$ is a contraction with semi-spectral measure $\E$ and $K$ is a trace class operator  on Hilbert space. Then
$$
\trace\Big(\,\,\iint_{\T\times\T}\big(\dg f\big)(\z,\t)\,d\E(\z)K\,d\E(\t)\Big)=
\int_\T f'(\z)\,d\mu(\z),
$$
where $\mu$ is a complex Borel measure on $\T$ defined by}
$
\mu(\D)=\trace(K\E(\D)).
$

\medskip

Recall that to prove this result, we used in \cite{MNP2} a representation of the divided difference $\dg f$ in the form \rf{Haagerpre} subject to the condition
\rf{proisumm}. We will prove in \S\;8 of this paper a similar result in the case of dissipative operators.

The following theorem was proved in \cite{MNP2}. Since we are going to use the construction given in the proof, we reproduce its proof here.

\begin{thm}
\label{fsdm}
Let $T_0$ and $T_1$ be contractions on Hilbert space such that $T_1-T_0\in\bS_1$.
Then there exists a complex Borel measure $\nu$ on $\T$ such that trace formula
{\em\rf{maksborme}}
%\bay
%\label{formsledlyaszha}
%\trace\big(f(T_1)-f(T_0)\big)=\int_\T f'(\z)\,d\nu(\z)
%\ey
holds for every $f$ in $\OLA$.
\end{thm}

\Pf Let $K\df T_1-T_0$. Consider the family of contractions $T_t\df T_0+tK$, $0\le t\le1$.
By Theorem 4.1 of \cite{MNP2} stated above, the function $t\mapsto f(T_t)$ is differentiable in the strong operator topology. Put
$$
Q_t\df\lim_{s\to0}\frac1s\big(f(T_{t+s})-f(T_t)\big)
=\iint_{\T\times\T}
\big(\dg f)(\z,\t)\,d\E_t(\z)K\,d\E_t(\t).
$$

Since $K\in\bS_1$ and $f\in\OLA$, we have
$$
Q_t\in\bS_1,\;0\le t\le1,\quad\mbox{and}\quad\sup_{t\in[0,1]}\|Q_t\|_{\bS_1}<\be.
$$

It follows from the definition of $Q_t$ that the function
$t\mapsto Q_tu$ is measurable for every vector $u$. Then 
the function $t\mapsto\trace (Q_tW)$ is measurable for an arbitrary bounded operator $W$. In other words the $\bS_1$-valued function $t\mapsto Q_t$ is
weakly measurable. Then it must be strongly measurable because $\bS_1$ is separable (see \cite{Y}, Ch. V, \S\:4) and
$$
f(T_1)-f(T_0)=\int_0^1Q_t\,dt,
$$
where the integral is understood in the sense of Bochner.

By Theorem 5.1 of \cite{MNP2} stated above, 
$$
\trace Q_t=\int_\T f'(\z)\,d\nu_t(\z),
$$
where $\nu_t$ is the complex Borel measure on $\T$ defined by
\bay
\label{nute}
\nu_t(\D)\df\trace\big(K\E_t(\D)\big),
\ey
for a Borel subset $\D$ of $\T$. Here $\E_t$ is the semi-spectral measure of $T_t$.

%It would be tempting to define the measure $\nu$ by
%$$
%\nu(\D)=\int_0^1\nu_t(\D)\,dt
%$$
%and conclude that formula \rf{maksborme} holds for this measure $\nu$.
%However, it is not quite clear why the function $t\mapsto\nu_t(\D)$ is measurable.
%
%To overcome this problem, we consider the dual space to the disk-algebra $\CA$. It can naturally be identified with the quotient space 
%$\mM/H^1_0$,
%where $\mM$ is the space of complex Borel measures on $\T$,
%$H^1_0=\{f\in H^1:~f(0)=0\}$ 
% and $H^1$ is the Hardy class.
%Consider the function $t\mapsto \dot\nu_t$, where $\dot\nu_t$ is the coset in 
%$\mM/H^1_0$ that corresponds to $\nu_t$.
%Let us show that the map $t\mapsto\dot\nu_t$ is continuous in the weak-star topology on $\mM/H^1_0$. Indeed, let $h\in\CA$. We have
%$
%\langle h,\dot\nu_t\rangle=\trace\big(Kh(T_t)\big).
%$
%The result follows from the fact that the map $t\mapsto h(T_t)$ is continuous in the operator norm. This can be proved easily by approximating $h$ by analytic polynomials. 
%
%Consider now the integral
%$
%\int_0^1\dot\nu_t\,dt.
%$
%It is an element of $\mM/H^1_0$, and so it is equal to $\dot\nu$ for some $\nu\in\mM$. It remains to observe that $\nu$ satisfies equality \rf{maksborme}. $\bl$
%
%\medskip

%{\bf Drugaya kontsovka.} 

Let $\mM$ be the space of complex Borel measures on $\T$. Clearly,
$$
\|\nu_t\|_\mM\le\|K\|_{\bS_1}.
$$
We consider the weak-star topology on $\mM$ induced by the space $C(\T)$ of continuous functions on $\T$. Let us show that the map $t\mapsto\nu_t$ is continuous in this topology. Indeed, let $h\in C(\T)$. It is easy to see that
$$
\langle h,\nu_t\rangle=\trace\left(K\int_\T h(\z)\,d\E_t(\z)\right),\quad h\in C(\T).
$$
We have to show that the map $t\mapsto \int_\T h\,d\E_t$ is continuous in the operator norm. Indeed, if $h=z^n$, then 
$$
\int_\T h(\z)\,d\E_t(\z)=\left\{\begin{array}{ll}
T_t^n,&n\ge0,\\[.2cm]
(T_t^*)^n,&n<0.
\end{array}\right.
$$
Thus the map $t\mapsto \int_\T h\,d\E_t$ is continuous for all trigonometric polynomials $h$.
It remains to approximate $h$ by trigonometric polynomials. 

We can define the complex Borel measure $\nu$ by $\nu\df\int_0^1\nu_t\,dt$. The integral can be understood as the integral of the function $t\mapsto\nu_t$ that is continuous in the weak-star topology. We have
$$
\trace\big(f(T_1)-f(T_0)\big)=\int_0^1\trace Q_t\,dt=
\int_0^1\left(\int_\T f'(\z)\,d\nu_t(\z)\right)\,dt
=\int_\T f'(\z)\,d\nu(\z)
$$ 
which completes the proof.
$\bl$

\begin{lem}
\label{anenut}
For every $t$ in $(0,1)$, the measure $\nu_t$ defined by {\em\rf{nute}} is absolutely continuous with respect to Lebesgue measure on $\T$.
\end{lem}

\Pf Let $t\in(0,1)$. If $T_t$ is a completely nonunitary contraction, then by the Sz.-Nagy theorem (see, \cite{SNF}, Th. 6.4 of Ch. 2), its minimal unitary dilation $U_t$ has spectral measure $E_t$ that is absolutely continious with respect to Lebesgue measure. Since the semi-spectral measure $\E_t$ is a compression of $E_t$, it is also absolutely continuous with respect to Lebesgue measure, and so the measure $\nu_t$ defined by \rf{nute} is absolutely continuous as well.

Suppose now that $T_t$ is not completely nonunitary. Let $\h$ be the Hilbert space on which the contractions $T_0$ and $T_1$ are defined
and let $\h=\h_\sharp\oplus\h_\flat$ be the decomposition of $\h$ such that $T_t\big|\h_\sharp$ is completely nonunitary while 
$T_t\big|\h_\flat$ is unitary. Consider the restrictions 
$T_0\big|\h_\flat$ and $T_1\big|\h_\flat$. We do not assume a priori that
$\h_\flat$ is an invariant subspace of $T_0$ and $T_1$, instead we consider these restrictions as operators from $\h_\flat$ to $\h$. Clearly, $T_t\big|\h_\flat$ is a convex combination of $T_0\big|\h_\flat$ and $T_1\big|\h_\flat$, and so for every unit vector $x$ in $\h_\flat$, $T_tx$ is a convex combination of $T_0x$ and $T_1x$. Since all vectors of norm 1 are extreme points of the unit ball of $\h$
and $\|T_tx\|=\|x\|=1$, it follows that $T_0x=T_1x$, and so 
$\h_\flat$ is indeed an invariant subspace of $T_0$ and $T_1$ and 
 $T_0\big|\h_\flat=T_1\big|\h_\flat=T_t\big|\h_\flat$.

Thus, $K\big|\h_\flat=\0$. It follows that 
$$
\nu_t(\D)=\trace\big(\E_t(\D)K\big)=\trace\big(\E^\sharp_t(\D)K\big),
$$
where $\E^\sharp_t$ is the semi-spectral measure of $T_t\big|\h_\sharp$, which is absolutely continuous by the Sz.-Nagy theorem.
This implies that the scalar measure $\nu_t$ is absolutely continuous. $\bl$

\begin{thm}
\label{funspesdvbszh}
Let $T_0$ and $T_1$ be contractions on Hilbert space. Then there exists a  function $\bs\xi$ in $L^1(\T)$ such that 
\bay
\label{fosleoLif}
\trace\big(f(T_1)-f(T_0)\big)=\int_\T f'(\z)\bs\xi(\z)\,d\z
\ey
for every $f$ in $\OLA$.
\end{thm}

A function $\bs\xi$ satisfying the conclusion of the theorem is called a
{\it spectral shift function} for the pair of contractions
$\{T_0,T_1\}$.

\medskip

\Pf 
Consider the measure $\nu_t$ defined by \rf{nute}.
By Theorem \ref{fsdm}, formula \rf{maksborme} holds with
$\nu=\int_0^1\nu_t\,dt$ and by Lemma \ref{anenut}, $\nu_t$ is absolute continuous
for $t\in(0,1)$. We have to prove that $\nu$ is absolutely continuous with respect to 
Lebesgue measure on $\T$.

Let $d\nu_t=g_t\,d\m$, $g_t\in L^1(\T)$. Then the $L^1$-valued function
$t\mapsto g_t$ is weakly measurable on $(0,1)$. Indeed, let $h\in L^\be(\T)$. Then the function
\bay
\label{tintgth}
t\mapsto\int_\T g_th\,dm
\ey
is measurable on $(0,1)$. Indeed, for continuous functions $h$, the function \rf{tintgth} is 
continuous. This is
exactly the continuity of the function $t\mapsto\nu_t$ in the weak-star topology which was established in the proof of Theorem \ref{fsdm}. In the general case we can approximate the $L^\be$ function $h$ by a uniformly bounded sequence of continuous functions that converges to $h$ almost everywhere. 

We can consider now the integral 
$
\int_0^1 g_t\,dt
$
which can be understood as a Bochner integral. Clearly, 
$d\nu=\Big(\int_0^1 g_t\,dt\Big)d\m$, and so $\nu$ is absolutely continuous.
$\bl$

\medskip

{\bf Remark 1.} If $f$ is a function in $\CA$ such that $f(T_1)-f(T_0)\in\bS_1$ whenever $T_0$ and $T_1$ are contractions with $T_1-T_0\in\bS_1$, then $f\in\OLA$. This follows from the corresponding fact for functions of unitary operators and the equality $\OLA=\OL(\T)\cap\CA$ (see \S\:\ref{OperLirazra}).

\medskip

{\bf Remark 2.} It is easy to see that the reasoning given in the proofs of Theorems
\ref{fsdm} and \ref{funspesdvbszh} lead to the proof of the following more general facts:

(i) Let $\{T_n\}_{n\ge1}$ be a sequence of contractions that converges to a contraction $T$ in the norm and let $K\in\bS_1$. Consider the complex Borel measures $\nu_n$ and $\nu$ on $\T$ defined by
$$
\nu_n(\D)\df\trace\big(\E_n(\D)K\big)\quad\mbox{and}\quad
\nu(\D)\df\trace\big(\E(\D)K\big),
$$
where $\E_n$ is the semi-spectral measure of $T_n$ and $\E$ is the semi-spectral measure of $T$. Then $\lim_{n\to\be}\nu_n=\nu$ in the weak-star topology on $\mM$.

(ii) Let $\{T_t\}_{t\in[0,1]}$ be a family of contractions that depends on $t$ continuously in the norm and let $K\in\bS_1$. Consider the complex Borel measures $\nu_t$,
$0\le t\le1$, defined by 
$$
\nu_t(\D)\df\trace\big(\E_t(\D)K\big),
$$
where $\E_t$ is the semi-spectral measure of $T_t$. Let $\nu$ be the measure defined by $\nu=\int_0^1\nu_t\,dt$, where the integral is understood as an integral of an $\mM$-valued  function continuous in the weak-star topology. Suppose that each $\nu_t$ is absolutely continuous with respect to Lebesgue measure and $d\nu_t=g_t\,d\m$,
$g_t\in L^1(\T)$, $t\in(0,1)$. Then $\nu$ is absolutely continuous and $d\nu=g\,d\m$, where 
$g=\int_0^1g_t\,dt$ and the last integral is a Bochner integral of the $L^1$-valued function $t\mapsto g_t$.

\medskip

Recall that it is easy to see that if $\bs\xi$ is a spectral shift function for a pair $\{T_0,T_1\}$ of contractions, then the set of all spectral shift functions
for $\{T_0,T_1\}$ can be parametrized by
$$
\bs\xi+h,\quad h\in H^1.
$$

The following fact allows us to see that if we consider the class of complex Borel measures $\nu$, for which trace formula \rf{maksborme} holds for all analytic polynomials $f$, we get the same class.

\begin{cor}
\label{merabsne}
Let  $\{T_0, T_1\}$ be  a pair of contractions in $\h$ with trace class difference.
Let  $\nu$ be a complex Borel measure, for which the trace formula 
\begin{equation*}
\trace(f(T_1) - f(T_0)) = \int_\T f'(\zeta)d\nu(\zeta)\,
\end{equation*}
holds for all analytic polynomials $f$. 
Then $\nu$ is absolutely continuous with respect to Lebesgue measure on $\T$.
\end{cor}

\Pf By Theorem \ref{funspesdvbszh}, there exists a function $\bs\xi$ in $L^1(\T)$ such that trace formula \rf{fosleoLif} holds for arbitrary analytic polynomials $f$.
Let $\l$ be the complex Borel measure defined by $d\l=d\nu-\bs\xi\,d\z$. Since
$$
\int_\T f'(\z)\,d\nu(\z)=\int_\T f'(\z)\bs\xi(\z)\,d\z
$$
for all analytic polynomials $f$,
it is easy to verify that the Fourier coefficients of $\l$ satisfy the equality $\widehat\l(j)=0$ for $j<0$, and so by the Brothers Riesz theorem (see \cite{Koo}, Sect. V.C.4), $\l$ is absolutely continuous. $\bl$

\medskip

Note that the first paper treated pairs $\{T_0, T_1\}$ of non-unitary (and
non-self-adjoint) operators is due to Langer \cite{La}. Assuming that $T_1- T_0
\in \bS_1$  and using the Riesz--Dunford functional calculus,  he proved the following
trace formula:
\bay
\label{Langer_for-la}
\trace\big(f(T_1) - f(T_0)\big)=  \frac{1}{2\pi i}\int_{\Gamma}
f'(\z)\log\Delta_{T_1/T_0}(z)\,dz
\ey
Here $\Delta_{T_1/T_0}$ is the perturbation determinant, $f$ is a holomorphic
function on a domain $\Omega$ containing $\sigma = \sigma(T_0) \cup \sigma(T_1)$, and
$\Gamma$ is a contour consisting of finitely many positively oriented, rectifiable,
simple closed curves contained in $\Omega$ which contain $\sigma$ in their
interiors.

The first  approach to trace formulae for the resolvents for a pair $\{T_0, T_1\}$ with
unitary $T_0$ and a contraction $T_1$ satisfying  $T_1- T_0 \in \bS_1$ and $D_{T_1} \df
(I- T_1^*T_1)^{1/2}\in \bS_1$ was given by  Rybkin \cite{Ryb87,Ryb89,Ryb94}. His
formulae involve a complex spectral shift function, which is A-integrable but not
Lebesgue integrable. We are going to improve his results on A-integrable sprctral shift below, see Theorem \ref{thm_A-Integ_trace}.

Later on  Adamyan and Neidhardt \cite{AN} proved that  there exists  a real spectral
shift function $\bs{\xi} \in L^1(\T)$ for a pair $\{T_0, T_1\}$ of contractions
satisfying $T_1- T_0 \in \bS_1$ and assuming in addition that  $I - |T_1| \in \bS^0_1$
and $I - |T^*_0| \in \bS^0_1$. Here  
$$
\bS^0_1 \df\Big\{A :~ \sum_j s_j(A)\big|\log
s_j(A)\big| < \infty\Big\} \subset \bS_1 
$$
and $\{s_j(A)\}_{j\ge0}$  is the sequence of
singular values of an  operator $A$. It is shown in [MNP2] that this result from
\cite{AN} is sharp, hence the only assumption $T_1- T_0 \in \bS_1$ does not ensure
existence of a real spectral shift for the pair $\{T_0, T_1\}$.

The existence of a complex  spectral shift $\bs{\xi} \in L^1(\T)$  for a pair $\{T_0,
T_1\}$ under an additional assumption can be deduced from the corresponding  result in [MN2] for pairs of maximal
dissipative operators. In full generality Theorem
\ref{funspesdvbszh} was proved in [MNP2] by using the technique of double operator integrals with Krein's perturbation
determinants and the technique of boundary triplets.

We proceed now to trace formula in terms of  A-integrals.

\medskip

{\bf Definition.}
A complex measurable  function $g$ on $\T$ is called A-{\it integrable} on $\T$
if
\bay
\label{ochen'malen'koe}
\lim_{t\to\be}t\,\m\{\z\in\T:~|g|>t\}=0
\ey
and the limit
$$
\lim_{t\to\be}\int\limits_{\{\z:|g(\z)|<t\}}g(\zeta)\,d\m
$$
exists.
In this case the limit is called the A-{\it integral} of $g$ and is denoted by
$$
({\rm A})\int_{\T}g(\zeta)d\m.
$$

It is well known that for every $g\in L^1(\T)$, its harmonic conjugate  $\widetilde g$
is always A-integrable but does not have to belong to $L^1(\T)$
(see \cite{Koo}).  

Consider now the Riesz projections of an $L^1$ function $g$:
$$
(\pp_+g)(\z)=\sum_{j\ge0}\widehat g(j)\z^j,\quad|\z|<1,\quad\mbox{and}\quad
(\pp_-g) =\sum_{j<0}\widehat g(j)\z^j,\quad|\z|>1.
$$ 
Both $\pp_+g$ and $\pp_-g$ have nontangential boundary values almost everywhere on $\T$ (see \cite{Koo}). We are going to use the same notation $\pp_+g$ and $\pp_-g$
for the corresponding boundary-value functions on $\T$.
It follows that the functions $\pp_+g$ and $\pp_-g$ are also A-integrable.
Note that the fact that $\widetilde g$  satisfies \rf{ochen'malen'koe} was proved by Kolmogorov, see \cite{Koo}.

Suppose that $\bs{\eta}$ is a function on $\T$ such that the functions $qf$
are A-integrable for arbitrary trigonometric polynomial $q$. We say that
$\bs{\eta}$ is a {\it generalized spectral shift function} for a pair 
of contractions $\{T_0, T_1\}$ with trace class difference if
$$
\trace\big(f(T_1)-f(T_0)\big)={\rm(A)}\int_\T f'(\z)\bs{\eta}(\z)\,d\z
$$
for every trigonometric polynomial $f$.

\begin{thm}
\label{thm_A-Integ_trace}
  Let $\{T_0, T_1\}$ be a pair of contractions on with trace class
  difference. Then

 {\em(i)} there exists a function $\bs\eta_-$ in $\pp_-L^1(\T)$
such that $\f\bs\eta_-$ is {\em A}-integrable on $\T$ for every $\f$ in $H^\be$ and
the trace formula
\bay
\label{A-Integ_trace_f-la}
\trace\big(f(T_1)-f(T_0)\big) ={\rm(A)}\int_{\T} f'(\z)\bs\eta_-(\z)\,d\z
\ey
holds for every $f$ in $\OLA$;

{\em(ii)} there exists a real function $\bs\xi_{\rm r}$ on $\Bbb T$ such $\f\bs\xi_{\rm r}$ is {\em A}-integrable on $\T$ for every $\f$ in $H^\be$ and the trace formula 
$$
\trace\big(f(T_1)-f(T_0)\big) ={\rm(A)}\int_{\T} f'(\z)\bs\xi_{\rm r}(\z)\,d\z
$$
holds for every  $f$ in $\OLA$.
\end{thm}

\Pf
(i)  Let $\bs\xi\in L^1(\T)$ be a spectral shift function for the pair $\{T_0, T_1\}$
and let $f\in \OL_A$. We set $\bs\xi_+ = \pp_+ \bs\xi$ and
$\bs\xi_- = \pp_- \bs\xi$.   Since $f\in \OL_A$, it follows that
$f'\in H^{\infty}$. 
%By  Kolmogorov's  theorem (see \cite{Koo}),
%$\bs\xi_+\in H^p(\Bbb D)$ for each $p\in(0,1)$. Combining both facts, yields  $f'
%\bs\xi_+\in H^p \subset N^+$ (Smirnov class).
By Ul'yanov's theorem \cite{Ulyan61}, the functions $\f\bs\xi_{\pm}$ A-integrable for any $\f$ in $H^\be$. Thus, both functions $f'\bs\xi_{\pm}$ are (A)-integrable and the function $f'\bs\xi_{+}$
satisfies the hypotheses of Alexandrov's Theorem \cite{Alex81}, and so
$$
{\rm(A)}\int_\T f'(\zeta)\bs\xi_+(\z)\,d\z= 0.
$$
Combining this equality with trace formula \rf{maksborme}, we obtain
  \begin{equation}
    \begin{split}
\trace\big(&f(T_1)-f(T_0)\big) = \int_{\T} f'(\z)\bs\xi(\z)\,d\z  \\
=&
{\rm(A)}\int_{\T}f'(\z){\bs\xi_+}(\z)\,d\z
+{\rm(A)}\int_{\T}f'(\z){\bs\xi_-}(\z)\,d\z=
{\rm(A)}\int_{\T}f'(\z){\bs\xi_-}(\z)\,d\z.
\end{split}
\end{equation}
It suffices to set $\bs\eta_-\df{\bs\xi_-}$.

(ii) 
%Without loss of generality we assume that $\widehat {\bs\xi_{I}}(0)=0$ where
%$\bs\xi_{I} = \im\bs\xi$.   Otherwise one can replace $\bs\xi$  by $\bs\xi - \widehat
%{\bs\xi_{I}}(0)$ which is also a spectral shift function for the pair $\{T_0, T_1\}$.
We have
$$
\bs\xi=\re\bs\xi+{\rm i}\im\bs\xi =\re\bs\xi+{\rm i}\pp_-(\im\bs\xi)+{\rm
i}\pp_+(\im\bs\xi).
$$
As mentioned above, the functions
 $\f\pp_{\pm}(\im\bs\xi)$ are A-integrable whenever $\f\in H^\be$. Put
\bay 
\label{ksifl}
\bs\xi_{\rm r}\df\re\bs\xi + {\rm i}\left(\pp_-(\im\bs\xi) - \ov{\pp_-(\im\bs\xi)}\right),
\ey
we obtain  a real A-integrable function. Moreover, since $\im\bs\xi\in L^1(\T)$, we have
$$
 \bs\xi-\bs\xi_{\rm r}= {\rm i}\big(\pp_+(\im\bs\xi) + \ov{\pp_-(\im\xi)}\,\big) \in
 \pp_+L^1(\T).
$$
Next, repeating the above reasoning, we find that
that the function $f'(\bs\xi -\bs\xi_{\rm r})$  satisfies the hypotheses of
Alexandrov's theorem \cite{Alex81}. Hence,
\begin{equation}
\label{ksifl_new}
{\rm(A)}\int_{\T}f'(\z)(\bs\xi(\z) -\bs\xi_{\rm r}(\z))\, d\z= 0.
\end{equation}

Combining trace formula \rf{maksborme} with \eqref{ksifl} and \eqref{ksifl_new}, we arrive at
\begin{align*}
\trace\big(f(T_1)-f(T_0)\big)&=\int_{\T} f'(\z)\bs\xi(\z)\,d\z\\[.2cm]
&=
{\rm(A)}\int_{\T}f'(\z) {\bs\xi_{\rm r}(\z)}\,d\z
+{\rm(A)}\int_{\T}f'(\z)(\bs\xi(\z)
-\bs\xi_{\rm r}(\z))\,d\z\\[.2cm]
&=
{\rm(A)}\int_{\T}f'(\z){\bs\xi_{\rm r}(\z)}\,d\z. \quad\bl
\end{align*}

\medskip

Note that the first trace formulas of the form \eqref{A-Integ_trace_f-la} for pairs
$\{T_0, T_1\}$ with unitary $T_0 = U_0$ and contractive $T_1$ and an A-integrable
generalized spectral shift function was established by Rybkin in
\cite{Ryb87}-\cite{Ryb94}. Our Theorem \ref{thm_A-Integ_trace}
complements Rybkin's results for pairs $\{U_0, T_1\}$.  In particular, the existence of an
A-integrable  real generalized spectral shift function was never established earlier even for
pairs $\{U_0, T_1\}$. In connection with the existence of an A-integrable generalized
spectral shift for pairs of the form $\{A, A-\ri V\}$, where $A$ is a self-adjoint operator and $V$ is a nonnegative trace class operator, we also mention the recent
publication \cite{MakScZ1}.

\

\section{\bf The trace formula for unitary operators}
\setcounter{equation}{0}
\label{Uni}

\

Let $U_0$ and $U_1$ be unitary operators on Hilbert space such that
$U_1-U_0\in\bS_1$. Recall that in \cite{AP+} it was proved that there exists a finite signed Borel measure $\mu$ on $\T$ such that 
\bay
\label{fosleunimer}
\trace\big(f(U_1)-f(U_0)\big)=\ri\int_\T\z f'(\z)\,d\mu(\z)
\ey
for every operator Lipschitz function $f$ on $\T$. This was obtained by considering the parametric family $U_t=e^{\ri tA}U_0$, $0\le t\le1$, differentiating the function
$t\mapsto f(U_t)$
in the strong operator topology, representing the trace of the derivative
as
\bay
\label{slepro}
\ri\int_\T\z f'(\z)\,d\mu_t(\z),
\ey
where $\mu_t$ is a finite signed Borel measure on $\T$ defined in terms of the spectral measure of $U_t$. To obtain trace formula \rf{fosleunimer}, it remains to take the integral of \rf{slepro} over $[0,1]$ and define $\mu$ as
$\int_0^1\mu_t\,dt$.

Though the method given in \cite{AP+} works for all operator Lipschitz functions, it does not give the absolute continuity of $\mu$. In \cite{AP+} Krein's theorem \cite{Kr2} was used to conclude that $\mu$ is absolutely continuous.

It turns out that the technique used for contractions in the previous section allows us to prove that $\mu$ is absolutely continuous without using the Krein result.

Recall that for an operator Lipschitz function $f$ on $\T$, the condition $U_1-U_0\in\bS_1$ for unitary operators $U_0$ and $U_1$ implies that $f(U_1)-f(U_0)\in\bS_1$ (see, e.g., \cite{AP} and \S\:3 of this paper).

\begin{thm}
\label{fspsunio}
Let $U_0$ and $U_1$ be unitary operators on Hilbert space such that
$U_1-U_0\in\bS_1$. Then there exists a real function $\bs\xi$ in 
$L^1(\T)$ such that
\bay
\label{spsduniop}
\trace\big(f(U_1)-f(U_0)\big)=\int_\T f'(\z)\bs\xi(\z)\,d\z
\ey
for every  $f$ in $\OL(\T)$.
\end{thm}

\Pf Let $\mu$ be a Borel measure such that trace formula \rf{fosleunimer} 
holds for all $f$ in $\OL(\T)$. In particular, \rf{fosleunimer} holds for all
$f$ in $\OLA$. Since $U_0$ and $U_1$ are contractions, we can apply Corollary  
\ref{merabsne} to find that $\mu$ is absolutely continuous. Let $d\mu=w\,d\m$.
We can define $\bs\xi$ by $\bs\xi=(2\pi)^{-1}w$. Clearly, \rf{spsduniop} holds.
$\bl$

%Clearly, $\bs\xi$ can be chosen so that $\widehat{\bs\xi}(0)=0$. Since
%$\trace\big(f(U_1)-f(U_0)\big)$ is real for all real trigonometric polynomials
%$f$, the right-hand side of \rf{spsduniop} must be real. It can be shown that this implies that $\bs\xi$ is a real function. $\bl$

%Since $U_0$ and $U_1$ are contractions, we can apply the result of the previous section. Thus, there exists a function $\bs{\xi}_\natural$ in $L^1(\T)$ such that 
%$$
%\trace\big(f(U_1)-f(U_0)\big)=\int_\T f'(\z)\bs\xi_\natural(\z)\,d\z
%$$
%for every $f$ in $\OLA$. 
%
%Let $\mu$ be the Borel measure such that trace formula \rf{fosleunimer} 
%holds for all $f$ in $\OL(\T)$. It follows that
%$$
%\hat\mu(j)=2\pi\ri\,\hat\xi_\natural(j-1),\quad j\le0.
%$$
%It follows from the brothers Riesz theorem (see \cite{Koo}, Sect. V.C.4) that $\mu$ is absolutely continuous and the result follows from \rf{fosleunimer} that there exists a function $\bs\xi$ in $L^1(\T)$ such that \rf{spsduniop} holds for $f\in\OL(\T)$.
%
%Without loss of generality we may assume that $\hat{\bs\xi}(0)=0$.
%In this case $\bs\xi$ must be real-valued. This fact is
%well-known and follows very easily from the fact that for real trigonometric polynomials $f$, the operator $f(U_0)-f(U_1)$ is self-adjoint, and so 
%$$
%\int_\T f'(\z)\bs\xi(\z)\,d\z=\trace\big(f(U_1)-f(U_0\big)\in\R.\quad\bl
%$$

\medskip

{\bf Remark.} The class $\OL(\T)$ is the maximal class for the applicability of the trace formula for functions of unitary operators. Indeed, if $f$ is a continuous function on $\T$ such that $f(U_1)-f(U_0)\in\bS_1$ whenever $U_0$ and $U_1$ are unitary operators with trace class difference, then $f$ must be in $\OL(\T)$,
see \S\:\ref{OperLirazra}.

\

\section{\bf Trace formulae for self-adjoint operators}
\setcounter{equation}{0}
\label{samosoop}

\

As we have already mentioned in the introduction, Krein's problem to describe the maximal class of functions $f$, for which trace formula \rf{foslsso} is applicable for arbitrary self-adjoint operators $A_0$ and $A_1$ with trace class difference, was solved in \cite{Pe6}. It was shown in \cite{Pe6} that the maximal class of functions in question coincides with the class of operator Lipschitz functions on $\R$.
The method of \cite{Pe6} is based on differentiating the function 
$t\mapsto f(A_0+t(A_1-A_0))-f(A_0)$ in the Hilbert--Schmidt norm and leads to trace formula \rf{foslemer} for a signed Borel measure $\nu$. In this section we deduce the absolute continuity of $\nu$ from the corresponding result for unitary operators (see \S\;\ref{Uni}), which was deduced in turn from the absolute continuity of spectral shift in the case of functions of contractions (see \S\;\ref{szhali}).

However, we first consider the case of {\it resolvent comparable pairs}
$\{A_0,A_1\}$ of self-adjoint operators, i.e., self-adjoint operators satisfying
the assumption
\bay
\label{oprezsra}
(A_1+{\rm i}I)^{-1}-(A_0+{\rm i}I)^{-1}\in\bS_1.
\ey

\medskip

{\bf Definiton.} We say that a continuous function $f$ is {\it resolvent operator Lipschitz} if
\bay
\label{rezoLi}
\|f(A)-f(B)\|\le\const\big\|(A+\ri I)^{-1}-(B+\ri I)^{-1}\big\|
\ey
for arbitrary bounded self-adjoint operators $A$ and $B$. We denote by
$\OLr$ the class of resolvent operator Lipschitz functions on $\R$.

\begin{lem}
\label{ogrneogr}
If $f\in\OLr$, then $f$ is bounded and inequality {\em\rf{rezoLi}} holds for arbitrary not necessarily bounded self-adjoint operators $A$ and $B$.
\end{lem}

\Pf Consider the Cayley transforms 
$$
U=(A-\ri I)(A+\ri I)^{-1}\quad\mbox{and}\quad
V=(B-\ri I)(B+\ri I)^{-1}
$$
of $A$ and $B$. It is well known and it is easy to verify that $U$ and $V$ are unitary operators and
\bay
\label{razun}
U-V=-2\ri\big((A+\ri I)^{-1}-(B+\ri I)^{-1}\big).
\ey 
It is also well known that a unitary operator $U$ is the Cayley transform of a not necessarily bounded self-adjoint operator if and only if 1 is not an eigenvalue of $U$, see e.g., \cite{BS0}. Finally, a self-adjoint operator is bounded if and only if 1 is not in the spectrum of its Cayley transform.

To prove that $f$ is bounded, put $A=tI$, $t\in\R$, and $B=\0$. Then inequality
\rf{rezoLi} implies that $|f(t)-f(0)|\le\const$. Consider the function $\f$ defined by
\bay
\label{pesafu}
\f(\z)=f\big({\rm i}(1+\z)(1-\z)^{-1}\big),\quad\z\in\T\setminus\{1\}.
\ey
It can easily be deduced from \rf{razun} that inequality \rf{rezoLi} is equivalent to the inequality
\bay
\label{opLineUV}
\|\f(U)-\f(V)\|\le\const\|U-V\|
\ey
for arbitrary unitary operators $U$ and $V$ such that $1\not\in\s(U)$ and $1\not\in\s(V)$. Then $\f$ extends to an operator Lipschitz function on $\T$
(see \S\:3.1 of \cite{AP}). In particular, this implies that inequality
\rf{opLineUV} holds whenever $U$ and $V$ are unitary operators such that
$1\not\in\s_{\rm p}(U)$ and $1\not\in\s_{\rm p}(V)$. This together with \rf{razun}
means that inequality \rf{rezoLi} holds for arbitrary not necessarily bounded self-adjoint operators $A$ and $B$. $\bl$

\begin{thm}
\label{ekvOLRT}
Let $f$ be a continuous function on $\R$. The following are equivalent:

{\em(i)} $f$ is resolvent operator Lipschitz;

{\em(ii)} if $A$ and $B$ are self-adjoint operators such that
\bay
\label{rezsraAB}
(A+{\rm i}I)^{-1}-(B+{\rm i}I)^{-1}\in\bS_1,
\ey
then $f(A)-f(B)\in\bS_1$;

{\em(iii)} if $A$ and $B$ are self-adjoint operators satisfying {\em\rf{rezsraAB}},
then 
$$
\|f(A)-f(B)\|_{\bS_1}
\le\const\big\|(A+{\rm i}I)^{-1}-(B+{\rm i}I)^{-1}\big\|_{\bS_1};
$$

{\em(iv)} the function $\f$ defined by {\em\rf{pesafu}} extends to an operator Lipschitz function on $\T$.
\end{thm}

\Pf Let us first prove that (i) implies (iii). It was shown in the proof of Lemma
\ref{ogrneogr} that the function $\f$ defined by \rf{pesafu} extends to
an operator Lipschitz function on $\T$, and so 
$$
\|\f(U)-\f(V)\|_{\bS_1}\le\const\|U-V\|_{\bS_1}
$$
for arbitrary unitary operators $U$ and $V$ (see \cite{AP}, Ch. 1, see also \S\:\ref{OperLirazra} of this paper). Again, using
\rf{razun}, we can conclude that (iii) holds.

Obviously, (iii) implies (ii). 

Let us prove that (ii) implies (i). Using Cayley transform, we see that (ii) is equivalent to the condition that
$\f(U)-\f(V)\in\bS_1$ whenever $U$ and $V$ are unitary operators such that
$U-V\in\bS_1$ and $1\not\in\s_{\rm p}(U)$ and $1\not\in\s_{\rm p}(V)$. As before, $\f$ is defined by \rf{pesafu}. 

It suffices to show that $\f$ extends to a continuous function on $\T$ and $\f(U)-\f(V)\in\bS_1$ whenever $U$ and $V$ are arbitrary unitary operators with $U-V\in\bS_1$. Indeed, we have already mentioned in \S\:\ref{OperLirazra} that this would imply that $\f\in\OL(\T)$.

Assume that $\f$ cannot be extended to a continuous function on $\T$.
Then there exist sequences $\{\z_n\}_{n\ge0}$ and 
$\{\t_n\}_{n\ge0}$ of points of $\T\setminus\{1\}$ such that 
$$
\lim_{n\to\be}\z_n=1,\quad\lim_{n\to\be}\t_n=1,
\quad\sum_{n\ge0}|\z_n-\t_n|<\be\quad\mbox{but}\quad
\inf_{n\ge0}|\f(\z_n)-\f(\t_n)|>0.
$$
Consider the unitary operators $U$ and $V$ on $\ell^2$ defined by
$$
U(x_0,x_1,x_2,\cdots)=(\z_0x_0,\z_1x_1,\z_2x_2,\cdots)
$$
and
$$
V(x_0,x_1,x_2,\cdots)=(\t_0x_0,\t_1x_1,\t_2x_2,\cdots).
$$
Clearly, $U$ and $V$ have no point spectrum at 1, $U-V\in\bS_1$ but
the operator $\f(U)-\f(V)$ is not even compact. 
This allows us to assume that $\f\in C(\T)$.

Suppose now that $U$ and $V$ are unitary operators on a Hilbert space $\h$ such that $U-V\in\bS_1$.
We have to show that $\f(U)-\f(V)\in\bS_1$. We denote by $P$ and $Q$ the orthogonal projections onto the eigenspaces $\Ker(U-I)$ and 
$\Ker(V-I)$. 
%Clearly, we may assume that both eigenspaces are infinite-dimensional. Otherwise, we can consider the orthogonal sums of both operators with the identity operator in infinitely dimensional Hilbert space.

Let
$U_\flat$ and $V_\flat$ be the restrictions of $U$ and $V$ to their invariant
subspaces $(\Ker(U-I))^\perp$ and $(\Ker(V-I))^\perp$. Clearly, $U_\flat$ and $V_\flat$ are unitary operators on $(\Ker(U-I))^\perp$ and $(\Ker(V-I))^\perp$ that do not have point spectrum at 1.

Suppose that $\{\z_k\}_{k\ge0}$ is a sequence of points in $\T\setminus\{1\}$
such that
\bay
\label{dvaneva}
\sum_k|\z_k-1|<\be\quad\mbox{and}\quad\sum_k|\f(\z_k)-\f(1)|<\be.
\ey
Let $\{u_k\}$, $0\le k<\dim\Ker(U-I)$, and $\{v_j\}$, 
$0\le j<\dim\Ker(V-I)$, be orthonormal bases in the subspaces
$\Ker(U-I)$ and $\Ker(V-I)$. Note that the subspaces can be infinite dimensional, finite-dimensional or even trivial. In the last case the orthonormal basis has no vectors.
Consider the unitary operators
$$
U_\sharp=\sum_k\z_k(\cdot,u_k)u_k\quad\mbox{and}\quad
V_\sharp=\sum_j\z_j(\cdot,v_j)v_j
$$
on these two subspaces. We can define now the unitary operators 
$U_\natural$ and $V_\natural$ on $\h$ by
$$
U_\natural\df U_\flat\oplus U_\sharp\quad\mbox{and}\quad
V_\natural\df V_\flat\oplus V_\sharp.
$$
Clearly, $1\not\in\s_{\rm p}(U_\natural)$ and $1\not\in\s_{\rm p}(V_\natural)$.
We have
$$
U_\natural-V_\natural=(U_\natural-U)-(V_\natural-V)+(U-V).
$$
It is easy to see that
$$
U_\natural-U=\0\oplus(U_\sharp-I)\in\bS_1\quad\mbox{and}\quad
V_\natural-V=\0\oplus(V_\sharp-I)\in\bS_1.
$$
Since $U-V\in\bS_1$, it follows that $U_\natural-V_\natural\in\bS_1$,
and so $\f(U_\natural)-\f(V_\natural)\in\bS_1$.
We have
$$
\f(U)-\f(V)=(\f(U)-\f(U_\natural))-(\f(V)-\f(V_\natural))+
(\f(U_\natural)-\f(V_\natural)).
$$
It follows easily from the second inequality in \rf{dvaneva} that
$\f(U)-\f(U_\natural)\in\bS_1$ and \lb$\f(V)-\f(V_\natural)\in\bS_1$. Thus, 
$\f(U)-\f(V)\in\bS_1$. This implies that $\f\in\OL(\T)$ (see \S\:\ref{OperLirazra}),
and so $f\in\OLr$.

The equivalence of (i) and (iv) is established in the proof of Lemma \ref{ogrneogr}.
This completes the proof. $\bl$

\begin{thm}
\label{samoresov}
Let $\{A_0,A_1\}$ be a resolvent comparable pair of self-adjoint operators.
Then there exists a real measurable function
$\bs\xi$ on $\R$ satisfying
\bay
\label{ksi1+t2-1}
\int_\R\frac{|\bs\xi(t)|}{1+t^2}\,dt<\be
\ey
and such that for an arbitrary function $f$ in $\OLr$, 
$f(A_1)-f(A_0)\in\bS_1$
and the following trace formula holds:
\bay
\label{formus}
\trace\big(f(A_1)-f(A_0)\big)=\int_\R f'(t)\bs{\xi}(t)\,dt.
\ey
\end{thm}

The function $\bs\xi$ is called the {\it spectral shift function for the pair 
$\{A_0,A_1\}$}.

\medskip

\Pf Let $U_0$ and $U_1$ be the Cayley transforms of $A_0$ and $A_1$ and
%It  is easy to verify that
%$U_0$ and $U_1$ are unitary operators and
%\bay
%%\label{razun}
%U_1-U_0=-2\ri\big((A_1+\ri I)^{-1}-(A_0+\ri I)^{-1}\big)\in\bS_1
%\ey
%by \rf{oprezsra}. 
%It is also well known that a unitary operator $U$ is the Cayley transform of a not necessarily self-adjoint operator if and only if 1 is not an eigenvalue of $U$, see e.g., \cite{BS0}.
let now $f$ be a function of class $\OLr$. Consider the function $\f$ on $\T$ be  defined by \rf{pesafu}. By Theorem \ref{ekvOLRT}, $\f$ extends to an operator Lipschitz function on $\T$. It is easy to verify that
$f(A_0)=\f(U_0)$ and $f(A_1)=\f(U_1)$.
Thus, 
$$
f(A_1)-f(A_0)=\f(U_1)-\f(U_0)\in\bS_1.
$$

By Theorem \ref{fspsunio}, there exists a function $\bs\xi_{\rm u}$ in $L^1(\T)$ (a spectral shift function for $\{U_0,U_1\}$) such that
$$
\trace\big(\f(U_1)-\f(U_0)\big)=\int_\T\f'(\z)\bs\xi_{\rm u}(\z)\,d\z.
$$

 Put now
$$
\bs\xi(t)=\bs\xi_{\rm u}\left(\frac{t-{\rm i}}{t+{\rm i}}\right),\quad t\in\R.
$$
It is easily seen that $\bs\xi$ satisfies \rf{ksi1+t2-1} and trace formula
\rf{formus} holds. $\bl$

\medskip

The following result shows that the class $\OLr$ is the maximal class of functions,
for which trace formula \rf{formus} holds for arbitrary resolvent comparable self-adjoint operators $A_0$ and $A_1$.

\begin{thm}
\label{maksrezsra} 
Suppose that $f$ is a continuous function on $\R$ such that
$f(A_1)-f(A_0)\in\bS_1$ whenever $A_0$ and $A_1$ are self-adjoint operators satisfying {\em\rf{oprezsra}}. Then $f\in\OLr$.
\end{thm}

\Pf Indeed, it follows from \rf{razun} that the function $\f$ defined by  
\rf{pesafu} has the property that $\f(U)-\f(V)\in\bS_1$ whenever $U$ and $V$ are unitary operators that have no point spectrum at 1 and such that $U-V\in\bS_1$. 
We have shown in the proof of Theorem \ref{ekvOLRT} that $\f$ extends to a function in $\OL(\T)$ which means (see again Theorem \ref{ekvOLRT}) that $\f\in\OLr$. $\bl$

\medskip

Let us now proceed to the case when $A_0$ and $A_1$ are self-adjoint operators with trace class difference. The following theorem is a combination of Krein's theorem
\cite{Kr} and the main result of \cite{Pe6}. We prove it without using Krein's theorem. Instead, we deduce it from Theorem \ref{samoresov}. Recall (see \S\:\ref{OperLirazra}) that if $f\in\OL(\R)$ and $A_0$ and $A_1$ are self-adjoint operators with $A_1-A_0\in\bS_1$, then $f(A_1)-f(A_0)\in\bS_1$.

\begin{thm}
\label{forssamoyara}
Let $A_0$ and $A_1$ be self-adjoint operators such that $A_1-A_0\in\bS_1$. Then 
there exists a real function $\bs\xi$ in $L^1(\R)$ such that 
$$
\trace\big(f(A_1)-f(A_0)\big)=\int_\R f'(t)\bs\xi(t)\,dt
$$
for every operator Lipschitz function $f$ on $\R$.
\end{thm}

\Pf It is well known and it is easy to see that the assumption $A_1-A_0\in\bS_1$ implies \rf{oprezsra}. Then it follows from Theorem \ref{samoresov} that there exists a function $\bs\xi$ satisfying \rf{ksi1+t2-1} such that trace formula \rf{formus} holds for smooth compactly supported functions $f$.

On the other hand, the construction in \cite{Pe6} yields the existence of a finite signed Borel measure $\nu$ on $\R$ such that 
$$
\trace\big(f(A_1)-f(A_0)\big)=\int_\R f'(t)\,d\nu(t)
$$
for arbitrary operator Lipschitz functions $f$. It follows that 
$$
\int_\R f'(t)\bs\xi(t)\,dt=\int_\R f'(t)\,d\nu(t)
$$
for arbitrary smooth compactly supported functions $f$. This easily implies that the measure $\nu$ is absolutely continuous and $d\nu(t)=\bs\xi\,dt$. Thus 
$\bs\xi\in L^1(\R)$ and trace formula \rf{formus} holds for arbitrary operator Lipschitz functions $f$. $\bl$

\

\section{\bf The case of dissipative operators; resolvent comparable perturbations}
\setcounter{equation}{0}
\label{rezsradis}

\

We consider in this section a trace formula for functions of maximal dissipative operators $L_0$ and $L_1$ such that
\bay
\label{rezsra}
(L_1+\ri I)^{-1}-(L_0+\ri I)^{-1}\in\bS_1.
\ey
Such a trace formula was obtained in \cite{MNP2}. We deduce this trace formula here from the corresponding result for functions of contractions, see \S\:\ref{szhali}.

Let us remind the following fact that can be deduced from Theorem 3.9.9 of \cite{AP} and from the corresponding fact for functions of unitary operators, see 
\S\:\ref{OperLirazra}. 

\medskip

{\it Let $f\in\CA$. The following are equivalent:}

(i) $f\in\OLA$;

(ii) {\it for arbitrary contractions $T$ and $R$
such that $T-R\in\bS_1$, the operator $f(T)-f(R)$ is also of trace class}; 

(iii) {\it the inequality
\bay
\label{yadszhiLiTR}
\|f(T)-f(R)\|_{\bS_1}\le\const\|T-R\|_{\bS_1}
\ey
holds for arbitrary contractions $T$ and $R$ with trace class difference}.

\medskip

\medskip

{\bf Definition.}
By analogy with the definition given in \S\:\ref{samosoop} we say that an $H^\be(\C_+)$ function $f$ continuous on $\R$ belongs to the {\it class 
$\OLAr$ of resolvent operator Lipschitz functions analytic in} $\C_+$ if
\bay
\label{resdisOL}
\|f(L)-f(M)\|\le\const\big\|(L+\ri I)^{-1}-(M+\ri I)^{-1}\big\|
\ey
for arbitrary maximal dissipative operators $L$ and $M$.

\medskip 

Actually, we do not have to assume a priori that the function $f$ is bounded.
We could assume that $f$ is a Lipschitz function on $\clos\C_+$ that is analytic in 
$\C_+$. Then it could be shown as in Lemma \ref{ogrneogr} that $f$ must be bounded.

Consider the function $\f$ defined by
\bay
\label{analpesafu}
\f(\z)=f\big({\rm i}(1+\z)(1-\z)^{-1}\big),\quad\z\in\dd\setminus\{1\}.
\ey

The following result is an analog of Theorem \ref{ekvOLRT}.

\begin{thm}
\label{ekvOLRT}
Let $f$ be a function in $H^\be(\C_+)$ that is continuous on $\R$. The following are equivalent:

{\em(i)} $f\in\OLAr$;

{\em(ii)} if $L$ and $M$ are maximal dissipative operators such that
\bay
\label{rezsraAB}
(L+{\rm i}I)^{-1}-(M+{\rm i}I)^{-1}\in\bS_1,
\ey
then $f(L)-f(M)\in\bS_1$;

{\em(iii)} if $L$ and $M$ are maximal dissipative operators {\em\rf{rezsraAB}},
then 
\bay
\label{yadLipdis}
\|f(L)-f(M)\|_{\bS_1}
\le\const\big\|(L+{\rm i}I)^{-1}-(M+{\rm i}I)^{-1}\big\|_{\bS_1};
\ey

{\em(iv)} the function $\f$ defined by {\em\rf{analpesafu}} extends to a function of class $\OLA$.
\end{thm}

\Pf Consider the Cayley transforms of $L$ and $M$ defined by
\bay
\label{PreKeLM}
T=(L-\ri I)(L+\ri I)^{-1}\quad\mbox{and}\quad
R=(M-\ri I)(M+\ri I)^{-1}
\ey
of $L$ and $M$. It is well known that $T$ and $R$ are contractions.
It is also well known that a contraction is the Cayley transform of a maximal dissipative operator if and only if 1 is not its eigenvalue.
It is easy to verify that
\bay
\label{raszhat}
T-R=-2\ri\big((L+\ri I)^{-1}-(M+\ri I)^{-1}\big).
\ey

Let us prove that (i) implies (iv). Let $\f$ be the function defined by  
\rf{analpesafu}. Then $f(L)=\f(T)$ and $f(M)=\f(R)$. Thus, it follows from 
\rf{raszhat} and \rf{resdisOL} that
$$
\|\f(T)-\f(R)\|\le\const\|T-R\|
$$
for arbitrary contractions $T$ and $R$ with no point spectrum at 1. In particular, this is true if $T$ and $R$ are unitary operators with no point spectrum at 1.
In the proof of Lemma \ref{ogrneogr} it has been shown that this implies that
$\f|\T$ extends to an operator Lipschitz function on $\T$. Then $\f\in\CA$ and since
$\OLA=\OL(\T)\cap\CA$, it follows that $\f\in\OLA$.

It is obvious that we can reverse the reasoning and establish that (iv) implies (i).

To show that (iv) implies (iii), we observe that inequality \rf{yadLipdis} 
for arbitrary maximal dissipative operators $L$ and $M$
is equivalent to inequality \rf{yadszhiLiTR} for arbitrary contractions $T$ and $R$ with no point spectrum at 1 which is certainly true.

The fact that (iii) implies (ii) is trivial.

It remains to establish the implication (ii)$\Rightarrow$(iv). Clearly, (ii) is equivalent to the property that $\f(T)-\f(R)\in\bS_1$ whenever $T$ and $R$ are contractions with no point spectrum at 1 such that $T-R\in\bS_1$. In particular, this is true if $T$ and $R$ are unitary operators with no point spectrum at 1. It was established in the proof of Theorem \ref{ekvOLRT} that in this case $\f|\T\in\OL(\T)$. Again it follows from the equality $\OLA=\OL(\T)\cap\CA$ that $\f\in\OLA$.
$\bl$

\begin{thm}
\label{rezOLdis}
Let $L_0$ and $L_1$ be maximal dissipative operators satisfying {\em\rf{rezsra}}. Then there exists a complex measurable function $\bs\xi$ 
on $\R$ such that 
\bay
\label{L21=x2-1}
\int_\R|\bs\xi(x)|(1+x^2)^{-1}\,dx<\be
\ey
and for an arbitrary function $f$ of class $\OLAr$, $f(L_1)-f(L_0)\in\bS_1$ and
the following trace formula holds:
\bay
\label{forslerezsradisop}
\trace\big(f(L_1)-f(L_0)\big)=\int_\R f'(t)\bs\xi(t)\,dt.
\ey
\end{thm}

%{\bf Remark.} As in the case of resolvent comparable self-adjoint operators, under the hypotheses of the theorem, $f(L_1)-f(L_0)\in\bS_1$. This will be seen in the proof.
%
%\medskip

\Pf Consider the Cayley transforms $T_0$ and $T_1$ of $L_0$ and $L_1$, see \rf{PreKeLM}).
%$$
%T_0=(L_0-\ri I)(L_0+\ri I)^{-1}\quad\mbox{and}\quad
%T_1=(L_1-\ri I)(L_1+\ri I)^{-1}
%$$
%of $L_0$ and $L_1$. It is well known that $T_0$ and $T_1$ are contractions.
%%It is also well known that a contraction is the Cayley transform of a maximal dissipative operator if and only if 1 is not its eigenvalue.
%It is easy to verify that
%\bay
%%\label{raszhat}
%T_1-T_0=-2\ri\big((L_1+\ri I)^{-1}-(L_0+\ri I)^{-1}\big)\in\bS_1.
%\ey

Let $f\in\OLAr$ and let $\f$ be the function defined by \rf{analpesafu}.
Then $\f\in\OLA$.
It is easy to see that
\bay
\label{fL1L2fT1T2}
f(L_1)-f(L_0)=\f(T_1)-\f(T_0)\in\bS_1
\ey
(see \S\:\ref{OperLirazra}).
It follows from formula \rf{raszhat} that $T_1-T_0\in\bS_1$.

By Theorem \ref{funspesdvbszh}, there exists a function $\bs\xi_{\rm c}$ in $L^1(\T)$ (a spectral shift function for the pair $\{T_0,T_1\}$) such that 
\bay
\label{sleszhafo}
\trace\big(\f(T_1)-\f(T_0)\big)=\int_\T\f'(\z)\bs\xi_{\rm c}(\z)\,d\z.
\ey
Put now
$$
\bs\xi(t)=\bs\xi_{\rm c}\left(\frac{t-{\rm i}}{t+{\rm i}}\right),\quad t\in\R.
$$
It is easily seen that $\bs\xi$ satisfies \rf{L21=x2-1}. Trace formula
\rf{forslerezsradisop} follows immediately from \rf{fL1L2fT1T2} and \rf{sleszhafo} 
$\bl$

\medskip

{\bf Remark.}
It is easy to see that $\OLAr$ is the maximal class of functions, for which the 
following property holds
$$
(L_1+\ri I)^{-1}-(L_0+\ri I)^{-1}\in\bS_1\quad\Longrightarrow\quad
f(L_1)-f(L_0)\in\bS_1
$$
for maximal dissipative operators $L_0$ and $L_1$. Indeed, self-adjoint operators are maximal dissipative, and so as we have mentioned in \S\:\ref{samosoop}, $f\in\OLr$. It remains to use the equality $\OLA=\OL(\T)\cap\CA$ (see \S\:\ref{OperLirazra}).

\medskip

Note that the first result on the existence of a spectral shift function for  pairs $\{L_0,L_1\}$ of
maximal accumulative resolvent comparable operators satisfying the additional assumption
$\rho(L_0)\cap \C_+\not = \varnothing$ was obtained in \cite{MN}. Besides, under
this assumption, using the Langer method (see formula \eqref{Langer_for-la}), the
authors proved  formula \eqref{forslerezsradisop} for a class of functions holomorphic on the
union of the spectra $\sigma(L_0)\cup \sigma(L_1)$. In full generality Theorem \ref{rezOLdis} was
proved in [MNP2] by developing and combining the Birman--Solomyak DOI approach with the Krein approach of perturbation determinants.

\

\section{\bf The case of dissipative operators; additive trace class perturbations}
\setcounter{equation}{0}
\label{disaddi}

\

In this section we obtain an analog of the Lifshits--Krein trace formula for additive trace class perturbations of maximal dissipative operators.

Let $L$ be a maximal (not necessarily bounded) dissipative operator in a Hilbert space $\h$.
If $f$ is an operator Lipschitz function in $\clos\C_+$, we can define the operator $f(L)$ as follows. We have
$$
f(\z)=\frac{f_\ri(\z)}{(\z+\ri)^{-1}},\quad\z\in\C_+,
\qquad
\mbox{where}\quad f_\ri(\z)\df\frac{f(\z)}{\z+\ri}.
$$
Since $f$ is operator Lipschitz, the function $f_\ri$ is continuous in $\clos\C_+$ and continuous at infinity (see \rf{pronabe}). The (possibly unbounded) operator $f(L)$ can be defined by
$$
f(L)\df(L+\ri I)f_\ri(L)
$$
(see \cite{SNF}, Ch. IV, \S\;1). It follows from Th. 1.1 of Ch. IV of \cite{SNF}  that 
$$
f(L)\supset f_\ri(L)(L+\ri I),
$$
and so $D(f(L))\supset D(L)$.

\begin{thm}
\label{DaKrpodi}
Suppose that $L_0$ and $L_1$ are maximal dissipative operators in a Hilbert space such that $D(L_1)=D(L_0)$ and $L_1-L_0\in\bS_1$. Let $f\in\OLA(\C_+)$. Then the operator $f(L_1)-f(L_0)$
on $D(L_0)$ is given by
$$
f(L_1)-f(L_0)=\iint_{\R\times\R}
\frac{f(t)-f(s)}{t-s}\,d\E_{L_1}(t)(L_1-L_0)\,d\E_{L_0}(s),
$$ 
and so $f(L_1)-f(L_0)\in\bS_1$.
\end{thm}

\Pf
We are going to use Lemma 6.4 of \cite{AP1}, which gives a construction of a sequence of functions $\o_n$ in $H^\be(\C_+)$ with the following properties:

{\rm (i)} $\lim\limits_{n\to\be}\o_n(z)=1$ for every $z\in\C_+$,

{\rm (ii)} $\|\o_n\|_{H^\be}=1$ for every $n$,

{\rm (iii)} $({\rm i}+z)\o_n\in H^\be$ for every $n$,

{\rm (iv)}
$\lim\limits_{n\to\be}\|({\rm i}+z)\o_n^\prime(z)\|_{H^\be}=0$.

Then for an arbitrary maximal dissipative operator $L$ the following equality holds
$$
f(L)u=\lim_{n\to\be}(f\o_n)(L)u,\quad u\in D(L).
$$
Indeed, it follows from Theorem 1.1 of Ch. IV of \cite{SNF}  that for $u\in D(L)$,
\bay
\label{fonL}
(f\o_n)(L)u=f(L)\o_n(L)u\to f(L)u\quad\mbox{as}\quad n\to\be.
\ey

On the other hand by Theorem 6.6 of \cite{AP1}, 
\begin{align*}
%\label{fLM}
\lim_{n\to\be}\big((f\o_n)(L_1)&-(f\o_n)(L_0)\big)\\[.2cm]
&=\lim_{n\to\be}
\iint_{\R^2}\frac{\o_n(x)f(x)-\o_n(y)f(y)}{x-y}\,d\E_{L_1}(x)(L_1-L_0)\,d\E_{L_0}(y)\\[.2cm]
&=\iint_{\R^2}\frac{f(x)-f(y)}{x-y}\,d\E_{L_1}(x)(L_1-L_0)\,d\E_{L_0}(y).
\end{align*}
which together with \rf{fonL} implies the result. $\bl$

\begin{thm}
\label{spsdmedi}
Let $L_1$ and $L_0$ be maximal dissipative operators with 
$L_1-L_0\in\bS_1$. Then there exists a complex Borel measure $\nu$ 
on $\R$ such that the following trace formula holds
\bay
\label{fodisyadprime}
\trace\big(f(L_1)-f(L_0)\big)=\int_\R f'(x)\,d\nu(x)
\ey
for every  $f$ in $\OLA(\C_+)$.
\end{thm}

\Pf The proof is similar to the proof of Theorem 4.1 of \cite{MNP2} stated in 
\S\;\ref{szhali}.
Let $K\df L_1-L_0$. Consider the family $L_t\df L_0+tK$, $0\le t\le1$,
of maximal dissipative operators.

By Theorem 3.9.6 stated in \S\:\ref{OperLirazra} of \cite{AP}, there exist sequences $\{\f_n\}_{n\ge1}$
and $\{\psi_n\}_{n\ge1}$ of class $\CA(\C_+)$ such that
$$
\left(\sup_{\z\in\C_+}\sum_{n\ge1}|\f_n(\z)|^2\right)
\left(\sup_{\z\in\C_+}\sum_{n\ge1}|\psi_n(\z)|^2\right)
=\|f\|_{\OLA(\C_+)}^2
$$
and
$$
\frac{f(z)-f(w)}{z-w}=\sum_{n\ge1}\f_n(z)\psi_n(w),\quad z,~w\in\C_+.
$$

As in the proof of Theorem 4.1 of \cite{MNP2}, we can prove that 
$$
\lim_{h\to0}\frac1h\big(f(L_{s+h})-f(L_s)\big)
=\iint_{\R\times\R}\frac{f(x)-f(y)}{x-y}\,d\E_s(x)K\,d\E_s(y)
\df Q_s
$$
in the strong operator topology, where $\E_s$ is the semi-spectral measure of $L_s$. 

Indeed, by Theorem \ref{DaKrpodi},
\begin{align*}
\frac1h\big(f(L_{s+h})-f(L_s)\big)&=
\iint_{\R\times\R}\frac{f(x)-f(y)}{x-y}\,d\E_{s+h}(x)K\,d\E_s(y)\\[.2cm]
&=\sum_{n\ge1}\f_n(L_{s+h})K\psi_n(L_s),
\end{align*}
while
$$
\iint_{\R\times\R}\frac{f(x)-f(y)}{x-y}\,d\E_s(x)K\,d\E_s(y)=
\sum_{n\ge1}\f_n(L_{s})K\psi_n(L_s).
$$
Then following the proof of Theorem 4.1 of \cite{MNP2}, we see that it suffices to show that 
$$
\lim_{h\to\0}\|\f_n(L_{s+h})-\f_n(L)\|=0.
$$
This follows from the fact that the function $\f_n$, being in the class 
$\CA(\C_+)$, is uniformly operator continuous, i.e.,
for each $\e>0$ there exists $\d>0$ such that 
$$
\|\f_n(L)-\f_n(M)\|<\e
$$
whenever $M$ is a maximal dissipative operator satisfying $\|L-M\|<\d$.
The latter is a consequence of Theorem 7.2 of \cite{AP1} and an analog of Theorem 8.1 of \cite{AP*} for maximal dissipative operators.

As in the proof of Theorem \ref{fsdm}, we can show that
$$
f(L_1)-f(L_0)=\int_0^1Q_s\,ds.
$$
The integral on the right can be understood in the sense of Bochner. Hence,
$$
\trace\big(f(L_1)-f(L_0)\big)=\int_0^1\trace Q_s\,ds.
$$
On the other hand, as in the proof of Theorem \ref{fsdm}, we can show that
$$
\trace Q_s=\int_\R f'(x)\,d\nu_s(x),
$$
where $\nu_s$ is a complex Borel measure on $\R$ defined by
$$
\nu_s(\D)=\trace(\E_s(\D)K).
$$

The rest of the proof is similar to the proof what has been done in the proof of Theorem \ref{fsdm}. Let $\mM(\R)$ be the space of complex Borel measures on $\R$.
We consider the function $s\mapsto\nu_s$. 

Let us show that it is continuous in the weak topology 
$\s\big(\mM(\R),C_0(\R)\big)$, where $C_0(\R)$ is the space of continuous functions
on $\R$ vanishing at infinity. Indeed, it is easy to see that
$$
\langle h,\nu_s\rangle=\trace\left(K\int_\R h\,d\E_s\right),\quad h\in C_0(\R),
$$
and it suffices to show that the function $s\mapsto\int_\R h\,d\E_s$ is continuous in the norm. If $h\in C_0(\R)\cap\CA(\C_+)$, then $\int_\R h\,d\E_s=h(L_s)$. The fact
that the function $s\mapsto h(L_s)$ is norm continuous is a consequence of Theorem 7.2 of \cite{AP1} and Theorem 8.1 of \cite{AP*}.

On the other hand, if $h=\bar h_\natural$, where 
$h_\natural\in C_0(\R)\cap\CA(\C_+)$, then 
$\int_\R h\,d\E_s=\big(h_\natural(L_s)\big)^*$, and so the function
$s\mapsto\int h\,d\E_s$ is norm continuous. 

The continuity of the map $s\mapsto\int_\R h\,d\E_s$ for an arbitrary function $h$ in 
$C_0(\R)$ follows from the obvious fact that the set
$$
\left\{h_1+\bar h_2:~\;h_1,~h_2\in C_0(\R)\cap\CA(\C_+)\right\}
$$
is dense in $C_0(\R)$.

We can define now the measure $\nu$ as the integral $\int_0^1\nu_s\,ds$ of the weak-star continuous function $s\mapsto\nu_s$. It is easy to verify that formula \rf{fodisyadprime} holds. $\bl$

%We consider the function $s\mapsto\dot\nu_s$, where $\dot\nu_s$ is the coset in 
%$\mM(\R)/H^1(\C_+)$
%that corresponds to the measure $\nu_s$, $\mM(\R)$ is the space of complex Borel measures on $\R$ and $H^1(\C_+)$ is the Hardy class of functions analytic in $\C_+$. Then this function is continuous in the weak-$*$ topology on $\mM(\R)/H^1(\C_+)$. Indeed, for $h\in\CA(\C_+)$,
%we have $\langle h,\dot\nu_t\rangle=\trace\big(Kh(L_t)\big)$. Continuity follows now from the fact that $h$ is uniformly operator continuous, see the proof of Theorem \ref{spsdmedi}. Now it suffices to select a measure $\nu$ in the coset $\int_0^1\dot\nu_t\,dt$. $\bl$

\medskip

Now we are ready to establish the main result of this section.

\begin{thm}
Let $L_0$ and $L_1$ be maximal dissipative operators with 
$L_1-L_0\in\bS_1$. Then there exists a complex function $\bs\xi$ 
in $L^1(\R)$ such that the following trace formula holds
\bay
\label{slefodio}
\trace\big(f(L_1)-f(L_0)\big)=\int_\R f'(x)\bs\xi(x)\,dx
\ey
for every  $f$ in $\OLA(\C_+)$.
\end{thm}

\Pf By Theorem \ref{rezOLdis}, there exists a measurable function $\bs\eta$ on 
$\R$ such that $\int_\R|\bs\eta(x)|(1+x^2)^{-1}\,dx<\be$ and the
trace formula  
$$
\trace\big(f(L_1)-f(L_0)\big)=\int_\R f'(x)\bs\eta(x)\,dx
$$
holds for all functions $f$ in $\OLAr$.
It follows from Theorem  \ref{spsdmedi} that
$$
\int_\R\frac1{(\l-x)^{2}}\bs\eta(x)\,dx
=\int_\R\frac1{(\l-x)^{2}}\,d\nu(x),\quad\im\l<0.
$$

Consider the Radon complex measure $\mu$ defined by $d\mu=d\nu-\bs\eta\,d\m$, where $\m$ stands for Lebesgue measure on $\R$. Then
$$
\int_\R\frac1{(\l-x)^{2}}\,d\nu(x)=0,\quad\im\l<0.
$$
It can easily be deduced from
the brothers Riesz theorem (see e.g., Lemma 3.7 of \cite{MN} for details) that
 $\nu$ is absolutely continuous with respect to 
Lebesgue measure, and so trace formula \rf{slefodio} holds
with $\bs\xi$ being the Radon--Nikodym density of $\nu$, i.e.,
 $d\nu=\bs\xi\,d\m$. $\bl$

\medskip

{\bf Remark.} The class $\OLA(\C_+)$ is the maximal class, for which trace formula 
\rf{slefodio} holds for arbitrary pairs of resolvent comparable maximal dissipative operators. Indeed, this can be deduced from Theorem \ref{maksrezsra} and the fact that self-adjoint operators are maximal dissipative operators.

\medskip

First generalizations  of formula \eqref{foslsso} to the case of pairs $\{L_0, L_1\}$
for a maximal accumulative (dissipative) operator $L_1$ and a self-adjoint operator $L_0$  were obtained
by Rybkin \cite{Ryb84}, \cite{Ryb94} and Krein \cite{Kr87}. For instance, Krein
treated a pair $\{L_0, L_1\}$  with   $L_1 = L_0 - {\rm i}V$, $V\ge 0$, and $V\in
\bS_1$, and proved in \cite{Kr87} an analog of formula \eqref{foslsso} with right-hand
side $\int_\R f'(t)\,d\nu(t)$ for  a complex Borel measure $\nu$ 
and for functions $f$ of class ${\mathcal W}_1^+(\R)$.  Here ${\mathcal W}_1^+(\R)$ is the class of
functions $f$ whose derivative is the Fourier transform of a complex measure supported
in $[0,\be)$.

Trace formula \eqref{slefodio} for pairs $\{L_0,L_1\}$ of maximal accumulative operators with
trace class difference
was proved for $f$ in ${\mathcal W}_1^+(\R)$ in \cite{MN}. It was extended in \cite{MNP2} to the
class  $\OLA(\C_+)$  by developing both the approach by
Birman and Solomyak based on double operator integrals and Krein's
perturbation determinants approach.

\

\section{\bf A construction of an intermediate contraction}
\setcounter{equation}{0}
\label{promzh}

\

The purpose of this section is to show that under a certain additional condition on  a pair $\{T_0, T_1\}$ of
contractions with $T_1 - T_0 \in \bS_1$, there exists a contraction $T$ such that 
$T-T_0\in\bS_1$, and the pairs $\{T, T_0\}$ and $\{T,T_1\}$ have spectral shift functions
$\bs\xi_0$ and $\bs\xi_1$ satisfying
$$
\im\bs\xi_0\ge\0\quad\mbox{and}\quad\im\bs\xi_1\ge\0.
$$
Clearly, the function $\bs\xi$ defined by $\bs\xi\df\bs\xi_0-\bs\xi_1$ is a spectral shift function for the initial pair $\{T_0, T_1\}$. To achieve this, we use another parametric family of contractions that connects $T_0$ and $T_1$. A similar parametric family was used in \cite{AP+} for pairs of unitary operators.

We start with the following lemma:

\begin{lem}
\label{UTT}
Let $T$ be a contraction and let $U$ be a unitary operator on Hilbert space. If $U-I \in\bS_1$, then the pair $\{T,UT\}$ has a \emph{real spectral shift function} $\bs\xi$.
\end{lem}

\Pf
Since $U-I \in \bS_1$, it is easy to see that there exists  a trace class self-adjoint operator $A$ such that
$U = e^{\ri A}$. 
%Indeed, according to  the spectral  representation
%\begin{equation*}
%V = \int^{\pi}_{-\pi}e^{i\lambda}dE(\lambda).
%\end{equation*}
%where $E(\cdot)$ is the spectral measure of $V$.   Setting  $A :=
%\int^{\pi}_{-\pi}\lambda dE(\lambda) (= A^*)$ one gets $V = e^{iA}$. Notice that
%$\sigma(A) \subseteq [-\pi,\pi]$. Since  $(\bS_1 \ni) I - V = -2i \sin(A/2)e^{iA/2}$, we
%obtain $\sin(A/2) \in \bS_1$. In turn, this implies   $A \in \bS_1$.

Put
\begin{equation*}
T_t \df e^{\ri tA}T, \quad t \in [0,1].
\end{equation*}
Then $T_0 = T$ and $T_1 = UT$. 

We are going to use a combination of the methods given
in the proof of Theorem 4.1 of \cite{AP+} 
and the proof of  Theorem \ref{fsdm} of this paper. 
It can be shown that there exists a finite signed Borel measure $\nu$  on $\T$
such that
for any $f \in \OLA$ the following trace formula holds
\bay
\label{fosUTT}
\trace(f(UT) - f(T)) = \int_\T\ri\,\z f'(\zeta)\,d\nu(\zeta).
\ey

Indeed, as in the proof of Theorem \ref{fsdm} of \cite{AP+}, we can define the real Borel
measure $\nu_t$ on $\T$ by
$$
\nu_t(\D) = \trace(A\E_t(\D))=\trace\big((\E_t(\D))^{1/2}A(\E_t(\D))^{1/2}\big)\in\R,
$$
where $\E_t$ is the semi-spectral measure of $T_t$.

As in the proof of Theorem \ref{fsdm} of this paper, we  have
$$
\langle h,\nu_t\rangle=\trace\left(A\int_\T h(\z)\,d\E_t(\z)\right),\quad h\in C(\T),
$$
and we can conclude that the function $t\mapsto\nu_t$ is a continuous function in the weak-star topology of $\mM$ (see also Remark 2 of 
\S\:\ref{szhali}). Consider the integral $\int_0^1\nu_t\,dt\df\nu$ as an integral of a continuous function in the weak-star topology
of $\mM$. As in the proof of Theorem \ref{fsdm}, we can differentiate
the function $t\mapsto T_t$, express the derivative $Q_t$ of this function as
$$
Q_t=\ri\int_\T\z f'(\z)\,d\nu_t(\z)
$$
and integrate over the interval $[0,1]$. This leads to formula \rf{fosUTT}.

%To complete the proof, we apply Corollary \ref{merabsne}  to conclude that the
%measure $\nu$ is absolutely  continuous. 

By Corollary \ref{merabsne}, the measure $\nu$ must be absolutely continuous.
Let $\eta$ be a real $L^1$ function such that $d\nu=\eta\,d\m$, where $\m$ is normalized Lebesgue measure on $\T$. 
Clearly, $2\pi\ri\z\,d\m(\z)=d\z$, and so
\bay
\label{nuzdz}
\ri\,\z\,d\nu(\z)=\ri\,\z\eta(\z)\,d\m(\z)=\frac1{2\pi}\eta(\z)\,d\z.
\ey
This together with \rf{fosUTT} implies that
$$
\trace(f(UT) - f(T)) = \frac1{2\pi}\int_\T f'(\z)\eta(\z)\,d\z.
$$
It remains to define $\bs\xi$ by $\bs\xi\df\frac1{2\pi}\eta$.
$\bl$

\begin{lem}
\label{XTT}
Let $T$ be a contraction and let $X$ be a nonnegative contraction, i.e.,
$0 \le X \le I$. If $I-X \in \bS_1$ and $X$ is invertible, then the pair $\{T,XT\}$ has a purely imaginary spectral shift function $\bs\xi$ satisfying $\im\bs\xi\ge 0$.
\end{lem}

\Pf Clearly, $X \ge \varepsilon I > 0$ for some  $\e>0$.
Hence,  $X = \int^{1}_{\varepsilon}\lambda
dE(\lambda)$,  where $E$ 
is the spectral measure of $X$. Since $I-X \in \bS_1$, it follows that
the spectral measure $E$ is discrete with point masses at $\l_k$
$$
I-X=
\int^{1}_{\e}(1-\l)\,dE(\l)=
\sum_k(1-\l_k)E(\{\l_k\})
$$
and $\sum_k(1-\l_k)\dim E(\{\l_k\})<\be$.

Then the operator $D$ defined by
$$
-D\df\log X = \int^{1}_\e\log\l\,dE(\l) = \sum_k\log\l_k\,E(\{\l_k\}) 
= \sum_k \frac{\log\l_k}{1-\l_k}(1-\l_k)E(\{\l_k\})
$$
is a bounded nonnegative operator and  $X = e^{-D}$. Moreover, it is easy to see that $D \in \bS_1$.

Put
$$
T_t\df e^{-tD}T, \qquad t \in [0,1].
$$
Note that $T_0 = T$ and $T_1 = XT$. By analogy with the proof of Lemma \ref{UTT}, we can conclude  that there exists
a positive Borel measure $\nu$  such that  the following trace formula holds
   \begin{equation*}
\trace(f(XT) - f(T)) =-\int_\T \z f'(\z) d\nu(\z)
  \end{equation*}
for every $f$ in $\OLA$.  

Indeed, we can define the positive Borel measure $\nu_t$ by
$$
   \label{SSF_int_repres_3}
\nu_t(\Delta) = \trace(D\E_t(\Delta)) =\trace(D^{1/2}\E_t(\Delta)D^{1/2})\ge0
$$
for a Borel subset $\D$ of $\T$.

As in the proof of Lemma \ref{UTT} (see Remark 2 in \S\:\ref{szhali}), we can conclude that the function
$t\mapsto\nu_t$, $t\in[0,1]$, is continuous in the weak-star topology of $\mM$ and
define the positive Borel measure $\nu$ by $\nu=\int_0^1\nu_t\,dt$.

 By Corollary  \ref{merabsne}, the measure $\nu$ is
absolutely continuous. Let $\eta$ be an $L^1$ function such that $d\nu=\bs\eta\,d\m$.
We have
$$
\trace(f(XT) - f(T))=-\int_\T \z f'(\z) d\nu(\z)=
\frac1{2\pi}\ri\int_\T f'(\z)\bs\eta(\z)\,d\z
$$
(see \rf{nuzdz}). It remains to put $\bs\xi\df\frac{\ri}{2\pi}\eta$. $\bl$

\begin{lem}
\label{T1T0T}
Let $T_1$ and $T_0$ be invertible contractions such that 
%$0\in\rho(T_1)\cap \rho(T_0)$. \footnote{vvesti oboznachenie $\rho(T)$} and
$T_1 - T_0 \in \bS_1$. Then there exists an invertible contraction $T$ satisfying the following conditions:

{\em(i)} $T - T_0 \in \bS_1$;
%and $0\in \rho(T)$. 

{\em(ii)} the pair $\{T_0,T\}$ has a spectral shift function $\bs\xi_0$ satisfying 
$\im\bs\xi_0\ge\0$;

{\em(iii)} the pair $\{T,T_1\}$ has a spectral shift function $\bs\xi_1$ satisfying 
$\im\bs\xi_1\le\0$.

{\em(iv)} the function $\bs\xi$ defined by
$$
\bs\xi\df\bs\xi_0+\bs\xi_1
$$
is a spectral shift function for the pair $\{T_0,T_1\}$. 
\end{lem}

\Pf Consider the operator $R = T_1T^{-1}_0$. Clearly, $I - R \in \bS_1$ and  $I -
|R|^2=I-R^*R=I-R^*+R^*(I-R)\in\bS_1$. Hence, 
$I-|R|=(I-|R|^2)(I+|R|)^{-1}\in\bS_1$. Let $R =
U|R|$ be a  polar decomposition of $R$. Since $R$ is invertible, the operator $U$ is
unitary. Moreover, the inclusion 
$I-R\in\bS_1$ implies that $I-U=I-R+U(|R|-I)\in\bS_1$. Since $I-|R|\in\bS_1$,  there is a trace class self-adjoint operator $C$ such that
\begin{equation}
\label{polar_decom}
|R| = e^{C}\quad \text{and}\quad  R = Ue^{C}.
\end{equation}
Put 
$$
C_-\df-CE_C((-\infty,0))\quad\mbox{and}\quad
C_+\df CE_C([0,\infty).
$$
Clearly, $C_\pm \ge 0$, $C_\pm \in \bS_1$ and
$C = C_+ - C_-$. Let
 \begin{equation}
 \label{factorization_T}
T_2\df e^{-C_-}T_0\quad\mbox{and}\quad
R_0\df T_2T^{-1}_0 = e^{-C_-}.
  \end{equation}
The inclusion  $C_-\in\bS_1$  implies that
 \begin{equation}
 \label{6.15_I-R0_trace-cl}
T_0 - T_2 = (I - R_0)T_0 = (I - e^{-C_-} )T_0  \in \bS_1.
  \end{equation}
Since  $R_0$ is a nonnegative contraction satisfying
\eqref{6.15_I-R0_trace-cl}, the pair  $\{T_0,T_2\} = \{T_0,R_0T_0\}$ satisfies the
hypotheses of  Lemma \ref{XTT}, and so it
has a purely imaginary spectral shift function $\bs{\eta}_0$ satisfying
$\im\bs{\eta}_0\ge\0$.

Let now $T\df UT_2$. Clearly, $T_2-T=(I-U)T_2\in\bS_1$ and since
$U$ is unitary, Lemma \ref{UTT} guaranties, that  the pair $\{T_2,T\}$ has a real
spectral shift function $\bs\digamma$. This allows us to conclude that 
$\bs\xi_0=\bs\digamma+\bs\eta_0$ is a spectral shift function for
the pair $\{T_0,T\}$. Clearly,
$\im\bs\xi_0=\im\bs\eta_0\ge\0$.

Next, consider the pair $\{T,T_1\}$. Put $R_1 = T T^{-1}_1$. It follows from
\eqref{factorization_T}, \eqref{polar_decom},  and the identity  $T_1= RT_0$,   that
\begin{equation}
R_1=TT^{-1}_1=Ue^{-C_-}T_0 (RT_0)^{-1}=Ue^{-C_-}T_0 T^{-1}_0 e^{-C}U^*
 =Ue^{-C_+}U^* = e^{-UC_+U^*}.
\end{equation}
Hence, $T = e^{-UC_+U^*}T_1$. It follows that $T_1-T=(I - e^{-UC_+U^*})T_1\in\bS_1$. By Lemma \ref{XTT},  the pair $\{T_1,T\}$  has a spectral shift function
with nonnegative imaginary part. Thus, the pair $\{T,T_1\}$ has 
a spectral shift function
$\bs\xi_1$ such that $\im\bs\xi_1\le\0$. Setting $\bs{\xi}\df\bs\xi_0 +
\bs\xi_1$ we obtain a spectral shift function of the pair $\{T_1,T_0\}$ with 
the required properties. $\bl$

\medskip

For a contraction $T$ on a Hilbert space $\h$, we need the construction of the Sch\"affer matrix dilation $U^{[T]}$ on the two-sided sequence space $\ell^2_\Z(\h)$
of $\h$-valued sequences, see \cite{SNF}, Ch. 1, \S\:5.  Here we identify $\h$ with the subspace of sequences
$\{v_n\}_{n\in\Z}$ such that $v_j=\0$ for $j\ne0$.
Such a dilation does not have to be minimal. However, the advantage of this dilation is that it allows us to consider unitary dilations of contractions on $\h$ on the same space $\ell^2_\Z(\h)$. Recall that for a contraction $T$ we use the notation $D_T$ for the defect operator:
$$
D_T\df(I-T^*T)^{1/2}.
$$
Here is the block matrix of $U^{[T]}$:
\bay
\label{maShe}
U^{[T]}=
\left(\begin{matrix}
\ddots&\ddots&\vdots&\vdots&\vdots&\vdots&\vdots&\vdots&\iddots\\
\cdots&\0&I&\0&\0&\0&\0&\0&\cdots\\
\cdots&\0&\0&I&\0&\0&\0&\0&\cdots\\
\cdots&\0&\0&\0&D_T&-T^*&\0&\0&\cdots\\
\cdots&\0&\0&\0&T&D_{T^*}&\0&\0&\cdots\\
\cdots&\0&\0&\0&\0&\0&I&\0&\cdots\\
\cdots&\0&\0&\0&\0&\0&\0&I&\cdots\\
\iddots&\vdots&\vdots&\vdots&\vdots&\vdots&\vdots&\ddots&\ddots
\end{matrix}\right).
\ey
Here the entry $T$ is at the $(0,0)$ position. In other words, the entries 
$U^{[T]}_{j,k}$ of $U^{[T]}$ are given by
$$
U^{[T]}_{0,0}=T,\quad U^{[T]}_{0,1}=D_{T^*},\quad U^{[T]}_{-1,0}=D_T,
\quad U^{[T]}_{-1,1}=-T^*,\quad 
U^{[T]}_{j,j+1}=I\quad\mbox{for}\quad j\ne0,~-1,
$$
while all the remaining entries are equal to $\0$.

\begin{thm}
\label{fspsudi}
Let $T$ and $Q$ be contractions on a Hilbert space. Suppose that 
\lb$U^{[T]} - U^{[Q]}\in \bS_1$ and let $\bs\xi$ be a spectral shift function for 
the pair $\big\{U^{[T]},U^{[Q]}\big\}$. Then $T-Q\in\bS_1$ and
$\bs\xi$ is a spectral shift function for the pair $\{T,Q\}$.
\end{thm}

\Pf 
%Since the matrices $U^{[T]}$ and $U^{[Q]}$ have only four non-trivial entries,  it follows that
%$$
%U^{[T]} - U^{[Q]}\in \bS_1.
%$$
Let $\bs\xi$ be a a spectral shift function for  the pair 
$\big\{U^{[T]},U^{[Q]}\big\}$, i.e., the identity
\begin{equation}
\label{fUTUQ}
\trace\big(f\big(U^{[Q]}\big)-f\big(U^{[T]}\big)\big)=\int_\T f'(\z)\bs\xi(\z)\,d\z
\end{equation}
holds for every $f$ in $\OL(\T)$. In particular, \rf{fUTUQ} holds for every $f$ in 
$\OLA$.
Clearly, both  $U^{[T]}$ and $U^{[Q]}$ are upper triangular
matrices with the only nonzero diagonal entries $T$ and $Q$, and so for $n\ge0$, both  $\big(U^{[T]}\big)^n$ and $\big(U^{[Q]}\big)^n$ are upper triangular
matrices with the only nonzero diagonal entries $T^n$ and $Q^n$. Thus,
$$
\trace\big(\big(U^{[Q]}\big)^n-\big(U^{[T]}\big)^n\big) = \trace(Q^n - T^n),\quad n\ge0,
$$ 
and so
\begin{equation}
\label{fUTfUQfTfQ}
\trace\big(f\big(U^{[Q]}\big)-f\big(U^{[T]}\big)\big)= \trace\big(f(Q)-f(T)\big)
\end{equation}
for every $f$ in $\OLA$. Combining  identities \rf{fUTUQ} and \rf{fUTfUQfTfQ}, we
conclude that $\bs{\xi}$ is a spectral shift function for the pair $\{T,Q\}$.
$\bl$

\medskip

{\bf Remark.} Clearly, under the hypotheses of Theorem \ref{fspsudi}, the pair 
$\{T,Q\}$ has a real spectral shift function. Corollary 8.4 of \cite{MNP2} shows that there are pairs of contractions $\{T,Q\}$ with $T-Q\in\bS_1$ that have no real spectral shift functions. Theorem \ref{fspsudi} implies that for such pairs
$U^{[T]}-U^{[Q]}\not\in\bS_1$.

\begin{lem}
\label{veshfspsd}
Let $T$ be a Fredholm contraction of zero index. Then there is an invertible
contraction $Q$ such that $Q-T\in\bS_1$ and the pair
$\{T,Q\}$ has a real spectral shift function.
\end{lem}

\Pf
Since $T$ is a contraction with zero Fredholm index, there exists  a partial isometry
$V$ with initial space $\Ker T$ and final space $\Ker T^*$. We set $Q\df T + V$. 
Then $V^*V=P_{\Ker T}$ and $VV^*=P_{\Ker T^*}$, the orthogonal projections onto 
$\Ker T$ and $\Ker T^*$.

Obviously, 
$Q$ is an invertible
contraction  and $T-Q=-V\in\bS_1$. We have
$$
Q^*Q = T^*T + T^*V + V^*T + V^*V = T^*T + V^*V
$$
which yields
$$
D^2_Q=I-Q^*Q=I-T^*T-V^*V=D^2_T-V^*V=D^2_T-P_{\Ker T}.
$$
It follows that
$$
D_Q^2x=\left\{\begin{array}{ll}D^2_Tx,&x\perp\Ker T,\\[.2cm]
\0,&x\in\Ker T,
\end{array}\right.
$$
and so
$$
D_Qx=\left\{\begin{array}{ll}D_Tx,&x\perp\Ker T,\\[.2cm]
\0,&x\in\Ker T.
\end{array}\right.
$$
%$$
%(D^2_T-P_{\Ker T})^2=D^2_T-D_TP_{\Ker T}-P_{\Ker T}D^2_T+P^2_{\Ker T}.
%$$
Hence, $D_Q=D_T-P_{\Ker T}=D_T-V^*V$. 
Similarly, $D_{Q^*}=(I-QQ^*)^{1/2}=D_{T^*}-VV^*$. It follows that
\begin{equation}
\label{DTDT*}
D_T-D_Q=V^*V\in\bS_1\quad\text{and}\quad D_{T^*}-D_{Q^*}=VV^*\in \bS_1.
\end{equation}
Let $U^{[T]}$ and $U^{[Q]}$ be the Sch\"affer matrix unitary dilations of $T$ and $Q$ on $\ell^2(\Bbb Z, \h)$. Since the matrices $U^{[T]}$ and $U^{[Q]}$ have only four non-trivial entries, it follows from \rf{DTDT*} that
$$
U^{[T]} - U^{[Q]}\in \bS_1.
$$
 By Theorem \ref{fspsunio},  the pair $\{U^{[T]},U^{[Q]}\}$ has a real spectral shift
function. The result follows now from Theorem \ref{fspsudi}. $\bl$
%\begin{equation}
%%\label{fUTUQ}
%\trace\big(f\big(U^{[T]}\big)-f\big(U^{[Q]}\big)\big)=\int_\T f'(\z)\bs\xi(\z)\,d\z
%\end{equation}
%holds for every $f$ in $\OL(\T)$. In particular, \rf{fUTUQ} holds for every $f$ in 
%$\OLA$.
%Clearly, both  $U^{[T]}$ and $U^{[Q]}$ are upper triangular
%matrices with the only nonzero diagonal entries $T$ and $Q$, and so for $n\ge0$, both  $\big(U^{[T]}\big)^n$ and $\big(U^{[Q]}\big)^n$ are upper triangular
%matrices with the only nonzero diagonal entries $T^n$ and $Q^n$. Thus,
%$$
%\trace\big(\big(U^{[T]}\big)^n-\big(U^{[Q]}\big)^n\big) = \trace(T^n - S^n),\quad n\ge0,
%$$ 
%and so
%\begin{equation}
%%\label{fUTfUQfTfQ}
%\trace\big(f\big(U^{[T]}\big)-f\big(U^{[Q]}\big)\big)= \trace\big(f(T)-f(Q)\big)
%\end{equation}
%for every $f$ in $\OLA$. Combining  identities \rf{fUTUQ} and \rf{fUTfUQfTfQ}, we
%conclude that $\bs{\xi}$ is a real spectral shift function for the pair $\{T,Q\}$.
%$\bl$

\begin{thm}
\label{T1T0xi1xi0}
Let $T_0$ be a Fredholm contraction with zero index and let $T_1$
be a contraction such that $T_1-T_0 \in \bS_1$. Then there exists a contraction $T$
satisfying the following conditions:

{\em(i)} $T - T_0 \in \bS_1$;
%and $T_1 - T \in \bS_1$ 

{\em(ii)} the pair $\{T_0,T\}$ has a spectral shift functions $\bs\xi_0$
satisfying $\im\bs\xi_0\ge\0$;

{\em(iii)} the pair $\{T,T_1\}$ has a spectral shift functions $\bs\xi_1$
satisfying $\im\bs\xi_1\le\0$;

{\em(iv)} the function $\bs\xi$ defined by
$$
\bs\xi=\bs{\xi}_0+\bs{\xi}_1
$$
is a spectral shift function for the pair $\{T_0,T_1\}$.
\end{thm}

\Pf
Since $T_1-T_0\in\bS_1$, the contraction $T_1$ is also a Fredholm operator with zero
index. By Lemma \ref{veshfspsd},  there exist invertible contractions $Q_1$ and $Q_0$ such that $Q_j - T_j \in \bS_1$, $j=0,1$. Moreover, Lemma \ref{veshfspsd}
yields the existence of real spectral shift functions $\bs\digamma_1$ and 
$\bs\digamma_0$ of the
pairs $\{Q_1,T_1\}$ and $\{T_0,Q_0\}$.

Next, since  $Q_1-Q_0\in\bS_1$, the pair $\{Q_0,Q_1\}$  satisfies the hypotheses of
Lemma \ref{T1T0T}, and so there exists a contraction $T$ such that $T-Q_0\in\bS_1$ and $T-Q_1\in\bS_1$ and such that the pairs $\{Q_0,T\}$ and $\{T,Q_1\}$ have
spectral shift functions $\bs\eta_0$ and $\bs\eta_1$ satisfying
$\im\bs\eta_0\ge\0$ and $\im\bs\eta_1\le\0$. It is easy to verify that the functions
$$
\bs\xi_0\df\bs\eta_0+\bs\digamma_0\quad\mbox{and}\quad
\bs\xi_1\df\bs\eta_1+\bs\digamma_1 
$$
are spectral shifts functions for the pairs $\{T_0,T\}$ and $\{T,T_1\}$. Note that
$\im\bs\xi_0\ge\0$ and $\im\bs\xi_1\le\0$. Hence
$$
\bs\xi\df\bs\xi_0+\bs\xi_1=\bs\eta_0+\bs\eta_1+\bs\digamma_0+\bs\digamma_1
$$
is a spectral shift function for the pair $\{T_1,T_0\}$.
$\bl$

\medskip

It would be interesting to find out whether the conclusion of Theorem \ref{T1T0xi1xi0} remains valid for arbitrary pairs $\{T_0,T_1\}$ of contractions with $T_1-T_0\in\bS_1$. 

\begin{thm}
\label{TUmnipo}
Let $U$ be a unitary operator and  let $T$ be a contraction such that $T-U\in\bS_1$.  Then there exists a spectral shift function $\bs{\xi}$ for the
pair $\{U,T\}$ such that $\im\bs{\xi}\ge\0$.
\end{thm}

\Pf Let us
first assume that $T$ is invertible.  Consider the polar decomposition $T=V|T|$ of $T$. Clearly, $V$ is a unitary operator. Then $I-|T|\in\bS_1$. Indeed, let
$K=T-U\in\bS_1$. We have
\begin{align*}
I-|T|&=(I+|T|)^{-1}(I-|T|^2)=
(I+|T|)^{-1}\big(I-(U^*+K^*)(U+K)\big)\\[.2cm]
&=-(I+|T|)^{-1}(U^*K+K^*U+K^*K)\in\bS_1.
\end{align*}
It follows that
$$
V-U=V(I-|T|)+T-U\in\bS_1.
$$
Setting $W=VU^{-1}$, we see that $W$ is a unitary operator and
$W-I=(V-U)U^{-1}\in\bS_1$. Thus,
the pair $\{U,V\} = \{U,WU\}$ satisfies the hypotheses of Lemma \ref{UTT}, and so it has a real spectral shift function $\bs\digamma$.

Next setting $X=V|T|V^*$, we consider the pair $\{V,T\}=\{V,V|T|\}=\{V,XV\}$.  Clearly,
$$
X-I=V(|T| - I)V^*\in \bS_1 \quad \text{and}\quad  0\le X \le I.
$$
Therefore, by Lemma \ref{XTT}, the pair $\{V,T\}$ has a purely  imaginary
spectral shift function $\bs\eta$ satisfying $\im\bs\eta\ge\0$. Put $\bs\xi\df
\bs\digamma+\bs\eta$. It is easy to see that $\bs\xi$ is a spectral shift function of the pair $\{U,T\}$ and $\im\bs\xi=\im\bs\eta\ge\0$.

Suppose now that $0\in\s(T)$. Since $T-U\in\bS_1$, it follows that $T$ is a
Fredholm operator of zero index.  Then instead of Lemma \ref{UTT}, we apply Lemma
\ref{veshfspsd}. It implies the existence of an invertible contraction $Q$ such that  $Q-T\in\bS_1$ and such that  the pair $\{Q,T\}$ has a real
spectral shift function $\bs\digamma$.

On the other hand, by what we have already proved, the pair $\{U,Q\}$ has  a
spectral shift function $\bs\chi$ such that $\im\bs\chi\ge\0$. Setting $\bs\xi\df
\bs\digamma+\bs\chi$, we see that $\bs\xi$ is a spectral shift function for the pair $\{U,T\}$ and $\im\bs\xi\ge\0$.
$\bl$

\medskip

The following analog of Corollary \ref{TUmnipo} for dissipative operators  can be deduced from Corollary \ref{TUmnipo} by passing to Cayley transform.

\begin{cor}
\label{LAsamo}
Let $A$ be a self-adjoint operator (not necessarily bounded) and let $L$ be a maximal dissipative operator such that
$$
(A+\ri I)^{-1}-(L+\ri I)^{-1}\in\bS_1.
$$
Then the pair $\{A,L\}$ has a spectral shift function $\bs\o$ in 
$L^1\big(\R,(1+t^2)^{-1}\big)$ such that $\im\bs\o\ge\0$.
\end{cor}

Clearly, under the hypotheses of Corollary \ref{TUmnipo}, the pair $\{T,U\}$ has a spectral shift function $\bs\xi$ satisfying $\im\bs\xi\le\0$. This fact (as well as Corollary \ref{TUmnipo}) was established
by a quite different method in \cite{MNP2}.
Similarly, under the hypotheses of Corollary \ref{LAsamo}, the pair $\{L,A\}$ has a spectral shift function $\bs\o$ satisfying $\im\bs\o\le\0$. This was proved in
\cite{MN} under an additional assumtion and in \cite{MNP2} in the general case.

\

\section{\bf Appendix}
\setcounter{equation}{0}

\

In \S\:\ref{szhali} we associated with a contraction $T$ and a trace class operator $K$ the complex Borel measure $\nu$ on $\T$ defined by 
$\nu(\D)\df\trace\big(\E(\D)K\big)$, where $\E$ is the semi-spectral measure of $T$. 
For a fixed trace class operator $K$, we proved, that $T\mapsto\nu$ is a continuous map from the set of contractions equipped with the norm topology to the space
of complex Borel measures equipped with the weak-star topology (see the proof of Theorem \ref{fsdm} and Remark 2 in \S\:\ref{szhali}).
A similar problem was considered in 
\S\:\ref{disaddi} (see the proof of Theorem \ref{spsdmedi}) for maximal dissipative operators.
 
In this section we give an alternative approach to such problems. The approach is based on Sch\"affer matrix dilations of contractions, see \rf{maShe}. First, we establish a useful result on the structure of such dilations. Perhaps, this result can be known to experts. However, we were not able to find a reference.

\begin{thm}
\label{strumaShe}
Let $T$ be a contraction on a Hilbert space $\h$ and let $U^{[T]}$ be the Sch\"affer matrix unitary dilation of $T$ on $\ell_{\Z}^2(\h)$ defined by {\em\rf{maShe}}.
Then $\ell_{\Z}^2(\h)$ can be represented as the orthogonal sum
$\K\oplus\K^{\perp}$ of reducing subspaces of $U^{[T]}$ such that
\bay
\label{prmiundi}
\K=\clos\spn\big\{\big(U^{[T]}\big)^n\h:~n\in\Z\big\}
\ey
is the subspace of minimal unitary dilation of $T$ and the restriction of $U^{[T]}$ to $\K^\perp$ is unitarily equivalent to a bilateral shift.
\end{thm}

By a bilateral shift we mean the operator ${\mathcal S}_\cL$ on $\ell^2_\Z(\cL)$ given by
$$
{\mathcal S}_\cL\{x_n\}_{n\in\Z}=\{x_{n-1}\}_{n\in\Z}.
$$
Here $\cL$ is a Hilbert space.

\medskip

\Pf Put $\mathcal D_T\df\clos\Range D_T$ and 
$\mathcal D_{T^*}\df\clos\Range D_{T^*}$ and put $U\df U^{[T]}$. We have
$$
U\{x_n\}_{n\in\Z}=
\{\cdots,x_{-2},x_{-1},D_Tx_0 -T^*x_1,\boxed{Tx_0 + D_{T^*}x_1},x_2,\cdots\},\!\!
\quad\{x_n\}_{n\in\Z}\in\ell^2_\Z(\h),
$$
(see \rf{maShe}). Here and in what follows the framed entry corresponds to the term indexed by $0$. It is easy to verify that
$$
U^*\{x_n\}_{n\in\Z}\!=\!
\{\cdots\!,x_{-2},\boxed{D_Tx_{-1}\!+T^*x_0},\!-Tx_{-1} + D_{T^*}x_0,x_1,x_2,\cdots\},
\!\!\quad\{x_n\}_{n\in\Z}\in\ell^2_\Z(\h).
$$

Let $\K$ be the subspace defined by \rf{prmiundi}. Then
$$
\K=\Big\{\{x_n\}_{n\in\Z}\in\ell^2_\Z(\h):~x_n\in\mathcal D_T,~n\le-1,~\;h_0\in\h,~\;
x_n\in\mathcal D_{T^*},~n\ge1\Big\},
$$
see \cite{SNF}, Ch. 1, Sect. 5.

It is easy to see that
$$
\K^\perp=\Big\{\{x_n\}_{n\in\Z}\in\ell^2_\Z(\h):~x_n\in\Ker D_T,~n\le-1,~\;x_0=\0,~\;
x_n\in\Ker D_{T^*},~n\ge1\Big\}.
$$
Let
$$
x=\{\cdots,x_{-2},x_{-1},\boxed{\0},x_1,x_2,\cdots\}\in\K^\perp.
$$
It is easy to see that
\bay
\label{prasd}
Ux=\{\cdots,x_{-2},x_{-1},-T^*x_1,\boxed{\0},x_2,x_3,\cdots\}
\ey
and
\bay
\label{lesd}
U^*x=\{\cdots,x_{-3},x_{-2},\boxed{\0},-Tx_{-1},x_1,x_2,\cdots\}.
\ey
It is easy to see that $-T$ maps isometrically $\Ker D_T$ onto $\Ker D_{T^*}$, while 
$-T^*$ maps isometrically $\Ker D_{T^*}$ onto $\Ker D_T$. Identities \rf{prasd}
and \rf{lesd} show that $U\big|\K^\perp$ is unitarily equivalent to the bilateral shift on $\ell^2_\Z(\Ker D_T)$ or $\ell^2_\Z(\Ker D_{T^*})$ (the map $-T$ allows us to identify $\Ker D_T$ and $\Ker D_{T^*}$). $\bl$

%Let $V$ the partial isometry 
%on $\h$ such that $V\big|\Ker D_T=-T\big|\Ker D_T$ and $V\big|(\Ker D_T)^\perp=\0$.
%It is easy to see that $V^*\big|\Ker D_{T^*}=-T^*\big|\Ker D_{T^*}$ and 
%$V^*\big|(\Ker D_{T^*})^\perp=\0$. We have

\begin{cor}
If $T$ is a completely nonunitary contraction, then the spectral measure of $U^{[T]}$ is mutually absolutely continuous with Lebesgue measure on $\T$.
\end{cor}

\Pf By the Sz.-Nagy theorem (see \cite{SNF}, Th. 6.4 of Ch. II), the spectral measure of the minimal unitary dilation of $T$ is mutually absolutely continuous with Lebesgue measure. By Theorem \ref{strumaShe}, $U^{[T]}$ is the orthogonal sum of the minimal unitary dilation and a bilateral shift whose spectral measure is certainly mutually absolutely continuous with Lebesgue measure. $\bl$

\begin{cor}
\label{tochspTUT}
Let $T$ be a contraction. Then $\s_{\rm p}(T)=\s_{\rm p}\big(U^{[T]}\big)$.
\end{cor}

\Pf By Proposition 6.1 of Ch. II of \cite{SNF}, the point spectrum of $T$ coincides with the point spectrum of its minimal unitary dilation. The result follow from Theorem \ref{strumaShe} because the bilateral shift has no eigenvalues. $\bl$

\begin{cor}
\label{skhopspeme}
Let $K$ be a trace class operator. Suppose that $\{T_n\}_{n\ge1}$ is a sequence of contractions on a Hilbert space $\h$ that converges to a contraction $T$ in the norm. 
%Let $\nu_n$ and $\nu$ be complex Borel measures on $\T$ defined by
%$$
%\nu_n(\D)\df\trace\big(\E_n(\D)K\big)\quad\mbox{and}\quad
%\nu(\D)\df\trace\big(\E(\D)K\big),\quad\D\subset\T,
%$$
%where $\E_n$ and $\E$ are the semi-spectral measures of $T_n$ and $T$. 
Then
for each arc $\I$ of $\T$ whose endpoints do not belong to $\s_{\rm p}(T)$,
\bay
\label{skhspeme}
\lim_{n\to\be}\E_n(\I)=\E(\I)
\ey
in the strong operator topology.
%{\em(ii)} for every $\f\in C(\T)$,
%$$
%\lim_{n\to\be}\int_\T\f(\z)\,d\E_n(\z)=\int_\T\f(\z)\,d\E(\z)
%$$
%in the norm;
%
%{\em(iii)} the sequence $\{\nu_n\}_{n\ge1}$ converges to $\nu$ in the weak-star topology on the space of Borel measures.
\end{cor}

\Pf Let $U_n$ and $U$ be the Sch\"affer matrix dilations of $T_n$ and $T$. It follows easily from \rf{maShe} that 
\bay
\label{UnU}
\lim_{n\to1}U_n=U
\ey
in the norm. 

%Indeed, $\lim_{n\to\be}T_n^*T_n=T^*T$ in the norm, and so 
%$\f(T_n^*T_n)$ converges to $\f(T^*T)$ for an arbitrary continuous function $\f$ on 
%$[0,1]$ (it can be approximated by polynomials). Thus $D_{T_n}\to D_T$ in the norm.
%The same as true for $D_{T_n}$ which implies \rf{UnU}.

Suppose now that $\I$ is an arc of $\T$ whose endpoints do not belong to 
$\s_{\rm p}(T)$. By Corollary \ref{tochspTUT}, the endpoints of $\I$ do not belong to
$\s_{\rm p}(U)$.

Let $E_n$ and $E$ be the spectral measures of $U_n$ and $U$. Then
$$
\lim_{n\to\be}E_n(\I)=E(\I)
$$
in the strong operator topology, see \cite{AG}, Sect. 78. To get \rf{skhspeme}, it suffices to consider the compressions of the spectral measures to $\h$. $\bl$

\begin{cor}
Let $L_n$ and $L$ be maximal dissipative operators in a Hilbert space $\h$ such that
\bay
\label{prerazre}
\lim_{n\to\be}\big\|(L_n+\ri I)^{-1}-(L+\ri I)^{-1}\big\|=0.
\ey
Let $\E_n$ and $\E$ be the semi-spectral measures of $L_n$ and $L$.
The following hold:

{\em(i)} if $\a,\,\b\not\in\R\setminus\s_{\rm p}(L)$, then
$$
\lim_{n\to\be}\E_n\big((\a,\b)\big)=\E\big((\a,\b)\big)
$$
in the strong operator topology;

{\em(ii)} if $\f\in C_0(\R)$, then
$$
\lim_{n\to\be}\int_\R\f(t)\,d\E_n(t)=\int_\R\f(t)\,d\E(t)
$$
in the norm;

{\em(iii)} if $\f$ is a bounded continuous function on $\R$, then
$$
\lim_{n\to\be}\int_\R\f(t)\,d\E_n(t)=\int_\R\f(t)\,d\E(t)
$$
in the strong operator topology.
\end{cor}

It is easy to see that \rf{prerazre} holds if $\lim_{n\to\be}\|L_n-L\|=0$.

\medskip

\Pf (i) Let 
$$
T_n=(L_n-\ri I)(L_n+\ri I)^{-1}\quad\mbox{and}\quad T=(L-\ri I)(L+\ri I)^{-1}
$$
be the Cayley transforms of $L_n$ and $L$. Then $T_n$ and $T$ are contractions and 
by \rf{raszhat}, $\lim_{n\to\be}\|T_n-T\|=0$. Let $U_n$ and $U$ be the Scah\"affer matrix dilations of $T_n$ and $T$ (see \rf{maShe}). Then their inverse Cayley transforms
$$
A_n\df\ri(I-U_n)(I+U_n)^{-1}\quad\mbox{and}\quad A\df\ri(I-U)(I+U)^{-1}
$$
are self-adjoint resolvent dilations of $L_n$ and $L$. As we have observed in the proof of Corollary \ref{skhopspeme}, $\lim_{n\to\be}\|U_n-U\|=0$, and so
$$
\lim_{n\to\be}\big\|(A_n+\ri I)^{-1}-(A+\ri I)^{-1}\big\|=0
$$
by \rf{razun}. By Theorem 8.24 of \cite{RS},
\bay
\label{spmesamsoop}
\lim_{n\to\be}E_n\big((\a,\b)\big)=E\big((\a,\b)\big),\quad\a,~\b\not\in\s_{\rm p}(A),
\ey
where $E_n$ and $E$ are the spectral measures of $A_n$ and $A$.
It follows from Corollary \ref{tochspTUT} that $\s_{\rm p}(A)=\s_{\rm p}(L)$.
To complete the proof, it suffices to the compress identity \rf{spmesamsoop} to $\h$.

\medskip

(ii) Since $\lim_{n\to\be}\big\|(A_n+\ri I)^{-1}-(A+\ri I)^{-1}\big\|=0$, the result follows from Theorem 8.20 of \cite{RS}.

\medskip

(iii) Again, since $\lim_{n\to\be}\big\|(A_n+\ri I)^{-1}-(A+\ri I)^{-1}\big\|=0$, the result follows from Theorem 8.24 of \cite{RS}. $\bl$

\

\medskip

\footnotesize
\noindent
\begin{tabular}{p{5cm}p{4.2cm}p{5cm}}
M.M. Malamud  & H. Neidhardt & V.V. Peller \\
People's Friendship University   & Institut f\"ur Angewandte & Department of Mathematics\\
of Russia (RUDN University)& Analysis und Stochastik&Michigan State University  \\
6 Miklukho-Maklaya St, Moscow, &Mohrenstr. 39& East Lansing, Michigan 48824\\
6117198, Russian Federation&D-10117 Berlin&USA\\
&Germany& and\\
&& People's Friendship University\\
 && of Russia (RUDN University)\\
 && 6 Miklukho-Maklaya St, Moscow,\\
&& 6117198, Russian Federation
\end{tabular}

\end{document}